# UNCONDITIONAL BASES OF INVARIANT SUBSPACES OF A CONTRACTION WITH FINITE DEFECTS

SERGUEI TREIL

ABSTRACT. The main result of the paper is that a system of invariant subspaces of a (completely non-unitary) Hilbert space contraction $T$ with finite defects (rank$(I - T^*T) < \infty$, rank$(I - TT^*) < \infty$) is an unconditional basis (Riesz basis) if and only if it is uniformly minimal.

Results of such type are quite well known: for a system of eigenspaces of a contraction with defects $1-1$ it is simply the famous Carleson interpolation theorem. For general invariant subspaces of operators with defects $1-1$ such theorem was proved by V. I. Vasyunin. Then partial results for the case of finite defects were obtained by the author.

The present paper solves the problem completely.

## CONTENTS

Part of the paper was written in when the author was visiting Mathematical Science Research Institute, Berkeley, CA. Research at MSRI is supported in part by NSF grant DMS-9022140.

## Notation

|  |  |  |
|---|---|---|
| $\mathbb{C}$ | — | complex plane |
| $\mathbb{D}$ | — | open unit disc, $\mathbb{D} \stackrel{\text{def}}{=} \{\xi \in \mathbb{C} : |\xi| < 1\}$; |
| $\mathbb{T}$ | — | unit circle, $\mathbb{T} \stackrel{\text{def}}{=} \partial\mathbb{D}$; |
| $m, |\,.\,|$ | — | the normalized ($m(\mathbb{T}) = 1$) Lebesgue measure on $\mathbb{T}$; |
| $\hat{f}(k)$ | — | $k$th Fourier coefficient of the function $f$ defined on the unit circle $\mathbb{T}$; $\hat{f}(k) \stackrel{\text{def}}{=} \int_{\mathbb{T}} f(z) z^{-k} dm(z)$; |
| $P_E$ | — | orthogonal projection onto a closed subspace $E$; |
| $\mathcal{L}\{\dots\}$ | — | linear span of the set $\dots$; |
| $\text{span}\{\dots\}$ | — | closed linear span of a set $\dots$; |
| $H^2, H^2_-,$ $H^\infty$ | — | Hardy classes |

$$H^2 \stackrel{\text{def}}{=} \left\{ f \in L^2 \;:\; \hat{f}(k) = 0 \text{ for } k < 0 \right\}, \qquad H^2_- \stackrel{\text{def}}{=} L^2 \ominus H^2$$

$$H^\infty \stackrel{\text{def}}{=} \left\{ f \in L^\infty \;:\; \hat{f}(k) = 0 \text{ for } k < 0 \right\}.$$

Hardy classes can be identified with the spaces of functions analytic in the unit disc $\mathbb{D}$, see [1].

|  |  |  |
|---|---|---|
| $L^\infty(E \to E_*)$ | — | the set of all bounded measurable functions on $\mathbb{T}$ with values in the space of operators from $E$ to $E_*$; |
| $L^2(E),$ $H^2(E)$ | — | vector Lebesgue and Hardy classes consisting of all functions with values in a (separable Hilbert) space $E$; |
| $H^\infty(E \to E_*)$ | — | the operator Hardy class $H^\infty$; |
| $\mathbb{P}_+, \mathbb{P}_-$ | — | orthogonal projections onto $H^2$ and $H^2_-$, respectively; |
| $k_\lambda$ | — | normalized ($\|k_\lambda\|_2 = 1$) reproducing kernel in $H^2$, $$k_\lambda(z) := \frac{(1 - \|\lambda\|^2)^{1/2}}{1 - \overline{\lambda} z}.$$ |

# 1. Introduction

The main result of the paper is that a system of invariant subspaces of a (completely non-unitary) Hilbert space contraction $T$ with finite rank defect operators ($\operatorname{rank}(I - T^*T) < \infty$, $\operatorname{rank}(I - TT^*) < \infty$) is an unconditional basis if and only if it is uniformly minimal. Let us recall some definition and a background, that explains why this problem is interesting and important.

We restrict ourselves to the case of a separable Hilbert space, although most of the definitions make a perfect sense in a Banach space as well.

Let $\mathcal{H}$ be a separable Hilbert space, and let $\mathcal{E}_n$ be closed subspaces of $\mathcal{H}$. Here and below, unless the converse is stated, by a subspace we always mean a *closed* linear subspace. Suppose that the subspaces $\mathcal{E}_n$ are linearly independent, that means that among all finite linear combinations $\sum f_n$, $f_n \in \mathcal{E}_n$, only trivial ones (all $f_n = 0$) equal 0. Assume also that the system of subspaces $\mathcal{E}_n$ is complete, i.e. $\mathcal{H} = \operatorname{span}\{\mathcal{E}_n : n = 1, 2, \dots\}$, where $\operatorname{span}\{\dots\}$ stands for the closed linear span of the set $\{\dots\}$. If the collection $\{\mathcal{E}_n : n = 1, 2, \dots\}$ is finite, then it is a standard linear algebra fact that the system $\{\mathcal{E}_n : n = 1, 2, \dots\}$ is a basis for $\mathcal{H}$, that means any vector $f \in \mathcal{H}$ admits a unique representation

$$f = \sum f_n, \qquad f_n \in \mathcal{E}_n. \tag{1.1}$$

If the system $\{\mathcal{E}_n : n = 1, 2, \dots\}$ is infinite, the things become more complicated. It is natural to define a basis as a system such that any vector $f$ admit a representation (1.1) where the series converges in the norm of $\mathcal{H}$. And that is how bases in infinite-dimensional spaces are defined. Clearly, to be a basis the system has to be linearly independent and complete. Unfortunately (or may be fortunately) for infinite systems it is not sufficient, and we shall show that below.

There is a simple description of infinite bases due to Banach. To formulate it we need some notation. Given a complete linearly independent system of subspaces $\mathcal{E}_n$, one can define (on a dense set $\mathcal{L}\{\mathcal{E}_n : n = 1, 2, \dots\}$ consisting of all finite linear combinations of $f_n \in \mathcal{E}_n$) projections $\mathcal{P}_n$

$$\mathcal{P}_n\Big(\sum f_k\Big) = \sum_{k=1}^{n} f_k, \qquad f_n \in \mathcal{E}_n. \tag{1.2}$$

**Theorem 1.1** (Banach Basis Theorem). *A complete linearly independent system of subspaces $\mathcal{E}_n$, $n = 1, 2, \dots$ is a basis if and only if all projections $\mathcal{P}_n$ defined above are uniformly bounded, $\|\mathcal{P}_n\| \leq C < \infty$.*

In the definition of a basis the ordering of subspaces $\mathcal{E}_n$ is essential, because the convergence of a series usually depends on the order. If we require that in the decomposition (1.1) the series converges unconditionally, i.e. independently of order, we obtain what is called an *unconditional basis*.

Unconditional bases admit a description similar to given by the Banach Basis Theorem above. For a finite subset $\sigma$ of $\mathbb{N}$, one can define on $\mathcal{L}\{\mathcal{E}_n \;:\; n = 1, 2, \ldots\}$ projections $\mathcal{P}_\sigma$

$$\mathcal{P}_\sigma\Big(\sum f_k\Big) = \sum_{k \in \sigma} f_k, \qquad f_n \in \mathcal{E}_n.$$

**Theorem 1.2.** *A complete linearly independent system of subspaces $\mathcal{E}_n$, for $n = 1, 2, \ldots$, is an unconditional basis if and only if all projections $\mathcal{P}_\sigma$ defined above are uniformly bounded, $\sup_\sigma \|\mathcal{P}_\sigma\| \le C < \infty$, where supremum is taken over all finite subsets $\sigma$ of $\mathbb{N}$.*

The above two theorems hold in any Banach space. In a Hilbert space we have a wonderful notion of an orthogonal basis. It is natural to consider bases that are orthogonal "up to isomorphism", i.e. that can be transformed into orthogonal basis by a bounded invertible operator. Such bases are called Riesz bases. An equivalent definition is that a Riesz basis is a basis that is orthogonal in an equivalent Hilbert norm.

It turns out that for a Hilbert space a system of subspaces is an unconditional basis if and only if it is a Riesz basis, [1, 3]. So in the sequel we will use both terms. We need a couple more of important definitions. With each Riesz basis one can associate a bounded invertible operator $\mathcal{J}$ that maps this basis into an orthogonal one. Such an operator (which is clearly not unique) is called an *orthogonalizer* of the system of subspaces.

It is convenient to "normalize" an orthogonalizer by assuming that its restriction on each subspace $\mathcal{E}_n$ is an isometry. Such "normalized" orthogonalizer is now unique up to a unitary multiplier on the left, and in what follows by orthogonalizer we always mean a normalized orthogonalizer.

If $\mathcal{J}$ is a (normalized) orthogonalizer for a Riesz basis $\{\mathcal{E}_n\}_{n=1}^\infty$, then the quantity $C(\{\mathcal{E}_n\}_{n=1}^\infty) = \|\mathcal{J}\| \cdot \|\mathcal{J}^{-1}\|$ could serve as measure on non-orthogonality of $\{\mathcal{E}_n\}_{n=1}^\infty$.

1.1. **A little more geometry.** We need couple more notions that generalize linear independence to the case of infinite system of subspaces. We say that a system $\{\mathcal{E}_n\}_{n=1}^\infty$ is *minimal* if all projections $\mathcal{P}^n$ defined on $\mathcal{L}\{\mathcal{E}_n \;:\; n = 1, 2, \ldots\}$

$$\mathcal{P}^n\big(\sum f_k\big) = f_n, \qquad f_k \in \mathcal{E}_k$$

are bounded, and *uniformly minimal* if all $\mathcal{P}^n$ are uniformly bounded $\sup_n \|\mathcal{P}^n\| =: 1/\delta < \infty$. The constant $\delta$ is called the *constant of uniform minimality* of the system $\{\mathcal{E}_n\}_{n=1}^\infty$.

One can easily see that the constant of uniform minimality $\delta$ admits very simple geometrical interpretation, namely

$$\delta = \inf_n \inf_{f \in \mathcal{E}_n, \|f\|=1} \mathrm{dist}\{f, \mathcal{L}\{\mathcal{E}_k \;:\; k \ne n\}\} \tag{1.3}$$

In other words, $\delta$ is the greatest lower bound for the sine of the angle between $\mathcal{E}_n$ and $\mathcal{L}\{\mathcal{E}_k : k \neq n\}$ }.

If a system $\{\mathcal{E}_n\}_{n=1}^\infty$ is minimal, one can define the so called *biorthogonal* or *dual* system $\{\mathcal{E}'_n\}_{n=1}^\infty$,

$$\mathcal{E}'_n = (\mathcal{P}^n)^* \mathcal{H}.$$

It is easy to see that a minimal system is uniformly minimal if and only if its dual is uniformly minimal.

It is also trivial that any minimal system is linearly independent, any uniformly minimal system is minimal, any basis is uniformly minimal, and any Riesz basis is a basis. In general there are great gaps between any two of the above properties. For example consider a system of one dimensional spaces $\mathcal{E}_n = \mathcal{L}\{z^n\}$, $n \in \mathbb{N}$ in the weighted space $L^2(w)$ where $w$ is a non-negative $L^1$ function (weight) on the unit circle $\mathbb{T}$. Then the system is:

1. Linearly independent iff $w \not\equiv 0$;
2. Minimal iff $\log w \in L^1$;
3. Uniformly minimal iff $1/w \in L^1$;
4. Basis iff $w$ satisfies Muckenhoupt $(A_2)$ condition, i.e. iff

$$\sup_I \left(\frac{1}{|I|}\int_I w\right) \cdot \left(\frac{1}{|I|}\int_I w^{-1}\right) < \infty;$$

5. Unconditional (Riesz) Basis iff $w \in L^\infty$, $1/w \in L^\infty$.

Let us explain that a little bit. The statement 1 is trivial. Let us prove 2. The condition $\log w \in L^1$ is equivalent to

$$\text{dist}_{L^2(w)}(1, \mathcal{L}\{z^n : n \geq 1\} > 0,$$

that means the projection $\mathcal{P}^0$, $\mathcal{P}^0\left(\sum_0^\infty \alpha_n z^n\right) = \alpha_0 \cdot 1$ is bounded. Since for the projection $\mathcal{P}^1$, $\mathcal{P}^1\left(\sum_0^\infty \alpha_n z^n\right) = \alpha_1 \cdot z$ we have $\mathcal{P}^1 f = z\mathcal{P}^0 \overline{z}(I - \mathcal{P}^0)f$, the projection $\mathcal{P}^1$ is bounded. Similarly all projections $\mathcal{P}^n$, $\mathcal{P}^n\left(\sum_0^\infty \alpha_k z^k\right) = \alpha_n z^n$ are bounded.

To prove 3 one need to notice that the system of positive exponents $\{z^n\}_0^\infty$ is uniformly minimal if and only if the system of all exponents $\{z^n\}_{-\infty}^\infty$ is. But for the last system the norms of all projections $\mathcal{P}^n$, $\mathcal{P}^n \sum \alpha_k z^k = \alpha_n z^n$ coincide, and it is an easy exercise on Cauchy–Swartz inequality to check that $\|\mathcal{P}^0\|^2 = \int_\mathbb{T} w \cdot \int_\mathbb{T} w^{-1}$.

By Theorem 1.1 the system $\{z^n\}_0^\infty$ is a basis if and only if $\sup_n \|\mathcal{P}_n\| < \infty$, where $\mathcal{P}_n\left(\sum \alpha_k z^k\right) = \sum_{k=0}^n \alpha_k z^k$. The last is equivalent to the boundedness of the Riesz projection $P_+$, $P_+ \sum \alpha_n z^n = \sum_0^\infty \alpha_n z^n$ in the weighted space $L^2(w)$. That is equivalent to the boundedness of the Hilbert transform $T$, $T = -iP_+ + i(I - P_+)$ and it is well known (see [5]) that $T$ is bounded if and only if $w$ satisfy the Muckenhoupt $(A_2)$ condition.

To prove 5 notice that if $\{z^n\}_0^\infty$ is a Riesz basis, the norm $\|\sum \alpha_k z^k\|_{L^2(w)}$ is equivalent to the norm $\left(\sum |\alpha_k|^2\right)^{1/2} = \|\sum \alpha_k z^k\|_{L^2}$. The last is possible if and only if $w, 1/w \in L^\infty$.

As one can see here the gap between Riesz Bases and uniformly minimal systems is as big as the gap between $L^\infty$ and $L^1$

However sometimes the gap does not exist. For example, let us consider a system $\{k_\lambda\}_{\lambda \in \sigma}$ of (normalized) reproducing kernels in the standard Hardy space $H^2$

$$k_\lambda(z) = \frac{(1-|\lambda|^2)^{1/2}}{1-\overline{\lambda}z},$$

where $\sigma$ is a countable subset of the unit disk $\mathbb{D}$. Note that $\|k_\lambda\|_2 = 1$. Then the corresponding system of one dimensional subspaces $\mathcal{E}_\lambda = \mathcal{L}\{k_\lambda\}$ is a Riesz basis (in its closed linear span, not in all $H^2$) if and only if it is uniformly minimal (as it was mentioned above the "only if" part is trivial).

Let as explain that in little more details. It is not a difficult exercise to compute the norm of the projection $\mathcal{P}^\lambda$ which is the projection onto $\mathcal{L}\{k_\lambda\}$ with kernel span$\{k_\mu : \mu \in \sigma \setminus \{\lambda\}\}$. Namely, cf [1, Lecture VI]

$$\|\mathcal{P}^\lambda\| = \left(\prod_{\mu \in \sigma \setminus \{\lambda\}} |b_\mu(\lambda)|\right)^{-1},$$

where $b_\lambda$ is a Blaschke factor, $b_\lambda(z) \stackrel{\text{def}}{=} (|\lambda|/\lambda)(\lambda - z)(1 - \overline{\lambda}z)^{-1}$. The condition $\sup_\lambda \|\mathcal{P}^\lambda\| < \infty$ is then nothing but the famous Carleson condition

$$\inf_{\lambda \in \sigma} \prod_{\mu \in \sigma \setminus \{\lambda\}} |b_\mu(\lambda)| =: \delta > 0,$$

that is necessary and sufficient condition for free interpolation in $H^\infty$ (and $\delta$ is exactly the constant of uniform minimality of the system $\{k_\lambda : \lambda \in \sigma\}$, see Section 1.1). Therefore given a finite subset $\sigma_1 \subset \sigma$ one can find a function $\varphi \in H^\infty$, $\|\varphi\|_\infty \leq C = C(\delta)$, such that $\varphi(\lambda) = 1$ for $\lambda \in \sigma_1$ and $\varphi(\lambda) = 0$ for $\lambda \in \sigma \setminus \sigma_1$. Then

$$\mathbb{P}_+ \overline{\varphi} k_\lambda = k_\lambda \quad \lambda \in \sigma_1 \quad \text{and} \quad \mathbb{P}_+ \overline{\varphi} k_\lambda = 0 \quad \lambda \in \sigma \setminus \sigma_1$$

Therefore for any $f \in \mathcal{L}\{k_\lambda : \lambda \in \sigma\}$ the projection $\mathcal{P}_{\sigma_1}$ onto span$\{k_\lambda : \lambda \in \sigma_1\}$ with kernel span$\{k_\lambda : \lambda \in \sigma \setminus \sigma_1\}$ is given by the formula

$$\mathcal{P}_{\sigma_1} f = \mathbb{P}_+ \overline{\varphi} f$$

and therefore $\|\mathcal{P}_{\sigma_1}\| \leq \|\varphi\|_\infty \leq C = C(\delta)$. So by Theorem 1.2 the system $\{k_\lambda : \lambda \in \sigma\}$ forms a Riesz basis in its closed linear span.

Let us note that it is possible to obtain directly from the above Carleson interpolation condition that the system $\{k_\lambda : \lambda \in \sigma\}$ is a Riesz basis and then get the result about free interpolation in $H^\infty$ using Sarason's theorem (a particular case of the Sz.-Nagy–Foiaş commutant lifting theorem). This scheme was realized in [1, Lectures VI, VII].

The connection between free interpolation in $H^\infty$ and the geometry (an unconditional basis property) of corresponding sequences of reproducing kernels in $H^2$ was first independently noticed by Nikolskii and Pavlov [6] and Katsnelson [2] in late 60's.

## 2. The main result

As everybody familiar with the Sz.-Nagy–Foiaş functional model for a contraction know the system of reproducing kernels $\{k_\lambda : \lambda \in \sigma\}$ has very simple operator–theoretic interpretation. It is simply a general form (up to unitary equivalence) of a system of eigenvectors of a completely non-unitary contraction with defect indices 1–1.

Let us remind that a contraction $T$ is called completely non-unitary if there exists no reducing subspace such that the restriction of $T$ on this subspace is a unitary operator. Since the theory of unitary operators is very well developed it seems very natural to "cut off" such unitary part of $T$ (if it exists) and consider only completely non-unitary contractions.

For a contraction $T$ one can consider the so called defect operators $D = (I-T^*T)^{1/2}$ and $D_* = (I-TT^*)^{1/2}$ and the corresponding defect indices $d = \operatorname{rank} D$, $d_* = \operatorname{rank} D_*$.

It seems natural to consider a general invariant subspaces of a contraction, and also not only contraction with defect indices 1–1, but with finite defect indices. It was proved by Vasyunin [10] that for an operator with defect indices 1–1 a complete system of invariant subspaces is a Riesz basis if and only if it is uniformly minimal. Then it was proved by the author that the same result holds for a system of eigenvectors of a contraction with finite defects [7, 8], and for a system of invariant subspaces of a contraction whose characteristic function is co-inner [9].

The main result of this paper is the following theorem.

**Theorem 2.1.** *Let $T$ be a completely non-unitary Hilbert space contraction with finite defect indices $d$, $d_*$, and $\{\mathcal{E}_n : n = 1, 2, \dots\}$ be a complete system of $T$-invariant subspaces. Then the system $\{\mathcal{E}_n : n = 1, 2, \dots\}$ is a Riesz basis if and only if it is uniformly minimal*

*Remark 2.2.* Let us note that in the paper we do not discuss the completeness at all. There are many results about when a system of eigenvectors, or of generalized eigenvectors of an operator is complete. We are interested in the question when a system is a basis in its own closed linear span. So in the statement of the above theorem we can omit the assumption that the system $\{\mathcal{E}_n : n = 1, 2, \dots\}$ is complete. In this case the conclusion will be "the system is a Riesz basis in its own linear span"

As it was already said this theorem for an operator with co-inner characteristic function was proved by the author in [9]. The presentation below mostly follows [9].

The main difference is the author now understands how to cope with outer part of the characteristic function. It turns out that one does not need anything fancy to deal with the outer part: it can be shown that by pretty trivial reasons the influence of the outer part can be eliminated.

Practically all new (comparing with [9]) techniques is contained in Sections 6, 7. However, presenting only this sections and saying that one can easily modify the construction in [9] to get the result, is a sure way to make an unreadable paper that can be understood besides the author by a couple of experts at best. So for the sake of readability and to make the paper self-contained the author needs to repeat the main lines of [9].

## 3. Idea of the proof of the main result.

The very general idea of the proof is quite simple. We need the following well know theorem (see [1, 3]) about Riesz bases.

**Theorem 3.1.** *A complete system of subspaces $\mathcal{E}_n$, $n = 1, 2, \ldots$ is a Riesz basis if and only if it is uniformly minimal and the following two "imbedding theorems*

1. $\sum_n \|P_{\mathcal{E}_n} f\|^2 \leq C \|f\|^2$
2. $\sum_n \|P_{\mathcal{E}'_n} f\|^2 \leq C \|f\|^2$

*hold for the system and its biorthogonal (dual) system $\mathcal{E}'_n$, $n = 1, 2, \ldots$.*

The main part of the proof is to show that if a system of invariant subspaces $\mathcal{E}_n$, $n = 1, 2, \ldots$ (of a completely non-unitary contraction $T$ with finite defects — I will omit that in the sequel) is uniformly minimal, then the imbedding theorem 1 from the theorem holds. If we do that we are done, because the dual system $\mathcal{E}'_n$, $n = 1, 2, \ldots$ is also uniformly minimal (trivially) and $\mathcal{E}'_n$ are clearly invariant subspaces of the adjoint operator $T^*$, which is of course a completely non-unitary contraction with finite defects.

According to the Sz.-Nagy–Foiaş functional model (see [4]) any (completely non-unitary) contraction in a separable Hilbert space admit the following representation (unitarily equivalent to the following operator):

Consider two Hilbert spaces $E$, $E_*$ and construct a space $\mathcal{H} = H^2(E) \oplus L^2(E_*)$. Let $S$ denotes the *shift* operator in $\mathcal{H}$, i.e. the multiplication by the independent variable $z$,

$$S \begin{pmatrix} f(z) \\ g(z) \end{pmatrix} = \begin{pmatrix} zf(z) \\ zg(z) \end{pmatrix}$$

and let $S^*$ be the adjoint of $S$

$$S^* \begin{pmatrix} f(z) \\ g(z) \end{pmatrix} = \begin{pmatrix} (f(z) - f(0))/z \\ \overline{z}g(z) \end{pmatrix}$$

Sz.-Nagy–Foiaş functional model theory says that any completely non-unitary contraction can be represented as a restriction of a backward shift $S^*$ onto its invariant subspace $\mathcal{K} \subset H^2(E) \oplus L^2(E_*)$, where $\dim E = d = \operatorname{rank}(1 - T^*T)$, $\dim E_* = d_* = \operatorname{rank}(1 - TT^*)$,

So we reduced our main theorem (Theorem 2.1) to the following one

**Theorem 3.2.** *Let $\mathcal{E}_n$, $n = 1, 2, \ldots$ be a system of $S^*$-invariant subspaces in $H^2(E) \oplus L^2(E_*)$, $\dim E < \infty$, $\dim E_* < \infty$. Suppose that the system $\mathcal{E}_n$, $n = 1, 2, \ldots$, is uniformly minimal, and let $\delta$ be the constant of uniform minimality. Then the imbedding theorem*

$$\sum_n \|P_{\mathcal{E}_n} f\|^2 \leq C \|f\|^2, \qquad \forall f \in H^2(E) \oplus L^2(E_*)$$

*holds, where $C = C(\delta, \dim E, \dim E_*)$.*

To prove the theorem we first need to remind the reader the description of $S^*$-invariant subspaces.

## 4. $S$- AND $S^*$-INVARIANT SUBSPACES OF $H^2(E) \oplus L^2(E_*)$

We need the following well known (see [1]) description of $S$-invariant subspaces in $H^2(E) \oplus L^2(E_*)$. Let us represent the space $H^2(E) \oplus L^2(E_*)$ as the set $\begin{pmatrix} H^2(E) \\ L^2(E_*) \end{pmatrix}$ consisting of all columns $\begin{pmatrix} f \\ g \end{pmatrix}$, $f \in H^2(E)$, $g \in L^2(E_*)$.

**Theorem 4.1.** *Let $\mathcal{M}$ be an $S$-invariant subspace of $\begin{pmatrix} H^2(E) \\ L^2(E_*) \end{pmatrix}$, $S\mathcal{M} \subset \mathcal{M}$. Then there exist a contractive analytic operator function $\Theta \in H^\infty(E_1 \to E)$, an operator function $\Delta \in L^\infty(E_1 \to E_*)$ and a function $P \in L^\infty(E_* \to E_*)$ whose values are orthogonal projections in $E_*$, such that $\Theta^*\Theta + \Delta^*\Delta = I$, $\operatorname{Range}\Delta(\xi) \perp \operatorname{Range}P(\xi)$ a.e. on $\mathbb{T}$ and*

$$\mathcal{M} = \mathcal{M}_{\Theta,\Delta,P} = \begin{pmatrix} \Theta \\ \Delta \end{pmatrix} H^2(E_1) \oplus \begin{pmatrix} 0 \\ PL^2(E_*) \end{pmatrix}$$

**Corollary 4.2.** *Any $S^*$ invariant subspace $\mathcal{K} \subset \begin{pmatrix} H^2(E) \\ L^2(E_*) \end{pmatrix}$ admits a representation*

$$\mathcal{K} = \mathcal{K}_{\Theta,\Delta,P} = \begin{pmatrix} H^2(E) \\ L^2(E_*) \end{pmatrix} \ominus \mathcal{M}_{\Theta,\Delta,P}$$

*where $\mathcal{M}_{\Theta,\Delta,P}$ and $\Theta$, $\Delta$, $P$ are from the previous theorem.*

*Remark 4.3.* Clearly if functions $\Theta$, $\Delta$, $P$ satisfy the assumptions of Theorem 4.1 then the corresponding subspace $\mathcal{M}_{\Theta,\Delta,P}$ is $S$-invariant. So, Theorem 4.1 and Corollary 4.2 give us a complete description of $S$-invariant and $S^*$-invariant subspaces respectively.

*Remark 4.4.* If $E_* = \{0\}$ Theorem 4.1 and Corollary 4.2 give a complete description of $S$ and $S^*$-invariant subspaces of $H^2(E)$. In this case $\Theta$ is an inner function, i.e. $\Theta(\xi)$ is an isometry a.e. on $\mathbb{T}$, $\Delta \equiv 0$, $P \equiv 0$ and

$$\mathcal{M} = \mathcal{M}_\theta = \Theta H^2(E_1), \qquad \mathcal{K} = \mathcal{K}_\Theta = H^2(E) \ominus \Theta H^2(E_1)$$

## 5. Imbedding theorems

For inner functions $\Theta_n$ the following results that were proved in [9].

**Theorem 5.1.** *Let $\theta_n$, $n \in \mathbb{N}$ be a sequence of (scalar valued) inner functions in $H^\infty$. Let $\mathcal{K}_{\theta_n} := H^2 \ominus \theta_n H^2$ and let $P_{\theta_n}$ denote the orthogonal projection in $H^2$ onto $K_{\theta_n}$. Then the imbedding theorem*

$$\sum_n \|P_{\theta_n} f\|^2 \le C \|f\|^2 \quad \forall f \in H^2 \tag{5.1}$$

*holds if and only if*

$$\sup_{\lambda \in \mathbb{D}} \sum_n (1 - |\theta_n(\lambda)|^2) = C < \infty \tag{5.2}$$

*Moreover the norm of the imbedding operator can be estimated above by a constant depending only on $C$.*

The following generalization of Theorem 5.1 to the case of matrix valued inner functions also can be found in [9].

**Theorem 5.2.** *Let $\Theta_n$, $n \in \mathbb{N}$ be a sequence of matrix valued inner functions ($\Theta_n(\xi)$ is an isometry a.e. on $\mathbb{T}$) in $H^\infty(E_n \to E)$. Let $\mathcal{K}_{\Theta_n} := H^2 \ominus \Theta_n H^2$ and let $P_{\theta_n}$ denote the orthogonal projection in $H^2$ onto $K_{\Theta_n}$. Then the imbedding theorem*

$$\sum_n \|P_{\Theta_n} f\|^2 \le C \|f\|^2 \quad \forall f \in H^2 \tag{5.3}$$

*holds if and only if*

$$\sup_{\lambda \in \mathbb{D}} \sum_n (\|\mathbf{e}\|^2 - \|\Theta_n(\lambda)^* \mathbf{e}\|^2) \le C' \|\mathbf{e}\|^2 \quad \forall \mathbf{e} \in E \tag{5.4}$$

*Moreover the norm of the imbedding operator can be estimated above by a constant depending only on $C$ and $\dim E$.*

If we assume that determinant of non-square matrix is 0 by definition, then the condition (5.4) immediately follows from

$$\sup_{\lambda \in \mathbb{D}} \sum_n (1 - |\det \Theta_n(\lambda)|^2) \le C'' < \infty \tag{5.5}$$

and moreover $C'$ from equation (5.4) is estimated above by $C''$ from (5.5). Clearly, if condition (5.5) holds there are only finitely many non-square matrices among $\Theta_n$.

Our goal is to show that if the system $\mathcal{K}_{\Theta_n,\Delta_n,P_n}$, $n \in \mathbb{N}$ (functions $\Theta_n$ are not necessarily inner) is uniformly minimal, then condition (5.5) is sufficient for the imbedding from Theorem 3.2 and then to prove that this condition takes place.

## 6. Some auxiliary results.

Let $\mathcal{K}_{\Theta,\Delta,P} = (\mathcal{M}_{\Theta,\Delta,P})^\perp$ be $S^*$-invariant subspace of $\begin{pmatrix} H^2(E) \\ L^2(E_*) \end{pmatrix}$, see Section 4, Corollary 4.2.

Let $k_\lambda$ denote the normalized reproducing kernel in $H^2$,

$$k_\lambda(z) := \frac{(1-|\lambda|^2)^{1/2}}{1-\bar{\lambda}z}$$

(note that $\|k_\lambda\|_{H^2} = 1$).

**Lemma 6.1.** *Let* $\begin{pmatrix} f \\ g \end{pmatrix} \in \begin{pmatrix} H^2(E) \\ L^2(E_*) \end{pmatrix}$. *Then*

$$P_{\mathcal{M}_{\Theta,\Delta,P}} \begin{pmatrix} f \\ g \end{pmatrix} = \begin{pmatrix} \Theta \\ \Delta \end{pmatrix} \mathbb{P}_+(\Theta^* f + \Delta^* g) + \begin{pmatrix} 0 \\ Pg \end{pmatrix}.$$

*In particular if $f = k_\lambda e$, $e \in E$, then*

$$P_{\mathcal{M}_{\Theta,\Delta,P}} \begin{pmatrix} k_\lambda e \\ 0 \end{pmatrix} = \begin{pmatrix} \Theta \\ \Delta \end{pmatrix} k_\lambda \cdot \Theta(\lambda)^* e.$$

*Proof.* Let $\mathcal{M} := \mathcal{M}_{\Theta,\Delta,P}$. Recall that $\mathcal{M}$ admits the decomposition

$$\mathcal{M} = \mathcal{M}_{\Theta,\Delta,P} = \begin{pmatrix} \Theta \\ \Delta \end{pmatrix} H^2(E_1) \oplus \begin{pmatrix} 0 \\ PL^2(E_*) \end{pmatrix} =: \mathcal{M}_1 \oplus \mathcal{M}_2,$$

see Theorem 4.1. Therefore $P_\mathcal{M} = P_{\mathcal{M}_1} + P_{\mathcal{M}_2}$.

It is easy to see that

$$P_{\mathcal{M}_2} \begin{pmatrix} f \\ g \end{pmatrix} = \begin{pmatrix} 0 \\ Pg \end{pmatrix}$$

To write down a formula for $P_{\mathcal{M}_1}$ consider an isometry $V : H^2(E_1) \to \begin{pmatrix} H^2(E) \\ L^2(E_*) \end{pmatrix}$,

$$Vf = \begin{pmatrix} \Theta f \\ \Delta f \end{pmatrix}, \qquad f \in H^2(E_1).$$

Clearly $\mathcal{M}_1 = \text{Range}V$, therefore $P_{\mathcal{M}_1} = VV^*$. Direct computations yield that

$$P_{\mathcal{M}_1} \begin{pmatrix} f \\ g \end{pmatrix} = VV^* \begin{pmatrix} f \\ g \end{pmatrix} == \begin{pmatrix} \Theta \\ \Delta \end{pmatrix} \mathbb{P}_+(\Theta^* f + \Delta^* g)$$

The statement about $k_\lambda e$ follows immediately from the trivial equality

$$\mathbb{P}_+ \Theta^* k_\lambda e = k_\lambda \Theta(\lambda)^* e$$

□

**Corollary 6.2.** *Let $f \in H^2(E)$. Then*

$$\operatorname{dist}\left\{ \begin{pmatrix} f \\ 0 \end{pmatrix}, \mathcal{K}_{\Theta,\Delta,P} \right\} = \|\mathbb{P}_+\Theta^* f\|$$

*Proof.* Applying Lemma 6.1 and using the fact that the operator $\begin{pmatrix} \Theta \\ \Delta \end{pmatrix}$ is an isometry we get

$$\operatorname{dist}\left\{ \begin{pmatrix} f \\ 0 \end{pmatrix}, \mathcal{K}_{\Theta,\Delta,P} \right\} = \left\| P_{\mathcal{M}_{\Theta,\Delta,P}} \begin{pmatrix} f \\ 0 \end{pmatrix} \right\| =$$

$$= \left\| \begin{pmatrix} \Theta \\ \Delta \end{pmatrix} \mathbb{P}_+\Theta^* f \right\| = \|\mathbb{P}_+\Theta^* f\|$$

□

Remind that $k_\lambda$ denote the normalized reproducing kernel in $H^2$,

$$k_\lambda(z) := \frac{(1-|\lambda|^2)^{1/2}}{1-\overline{\lambda}z}$$

(note that $\|k_\lambda\|_{H^2} = 1$).

**Corollary 6.3.** *Let $e \in E$. Then*

$$\operatorname{dist}\left\{ \begin{pmatrix} k_\lambda e \\ 0 \end{pmatrix}, K_{\Theta,\Delta,P} \right\} = \|\Theta(\lambda)^* e\|$$

*Proof.* Trivial corollary of Lemma 6.1. □

In the following lemma we assume that for a non-square matrix function

$$\Theta \in H^\infty(E_1 \to E)$$

the determinant is defined by $\det \Theta \equiv 0$ if $\dim E_1 < \dim E$ and $\det \Theta \equiv 1$ if $\dim E_1 > \dim E$.

**Lemma 6.4.** *Suppose a system of subspaces $\mathcal{M}_{\Theta_n,\Delta_n,P_n}$, $n \in \mathbb{N}$ in $\begin{pmatrix} H^2(E) \\ L^2(E_*) \end{pmatrix}$, $\Theta_n \in H^\infty(E_n \to E)$ is uniformly minimal with the constant of uniform minimality at least $\delta$. Then there exists a positive constant $\varepsilon = \varepsilon(\delta, \dim E)$ such that any point $\lambda \in \mathbb{D}$ is covered by at most $d := \dim E$ sets $\{z \in \mathbb{D} : |\det \Theta_n(z)| < \varepsilon^d\}$. In particular, at most $d$ matrices $\Theta_n$ satisfy $\dim E_n < \dim E$.*

*Proof.* Given a point $\lambda \in \mathbb{D}$ let us suppose there exist more that $d := \dim E$ functions $\Theta_n$ such that $|\det \Theta_n(\lambda)| < \varepsilon^d$. Without loss of generality one can assume that these functions are $\Theta_1, \Theta_2, \ldots, \Theta_{d+1}$ (part of these functions can be non-square).

For $k = 1, 2, \ldots, \dim E + 1$, one can find vectors $e_k \in E$ such that $\|\Theta_k(\lambda)^* e_k\| \leq \varepsilon$, $k = 1, 2, \ldots, \dim E + 1$. Note that if a function $\Theta_n$ is non-square and $\dim E_n < \dim E$, one can find a vector $e_n$ such that $\Theta_n(\lambda)^* e_n = 0$. By Corollary 6.3,

$$\mathrm{dist}\left\{\begin{pmatrix} k_\lambda e_k \\ 0 \end{pmatrix}, \mathcal{K}_{\Theta_k, \Delta_k, P_k}\right\} < \varepsilon.$$

Subspaces $\mathcal{K}_{\Theta_k, \Delta_k, P_k}$ are uniformly minimal, so for small enough $\varepsilon = \varepsilon(\delta, d)$ vectors $k_\lambda e_k$, $k = 1, 2, \ldots, \dim E + 1$ are uniformly minimal as well, and so are the vectors $e_k$. In particular vectors $e_k$ are linearly independent. But is is impossible because we have more that $\dim E$ vectors. □

**Lemma 6.5.** *Let $f \in L^2(E_*)$. Then*

$$\mathrm{dist}\left\{\begin{pmatrix} 0 \\ f \end{pmatrix}, \mathcal{K}_{\Theta, \Delta, P}\right\} = (\|Pf\|^2 + \|\mathbb{P}_+ \Delta^* f\|^2)^{1/2}$$

*Proof.* Trivial Corollary of Lemma 6.1. □

The following lemma can be viewed as a "boundary analogue" of Lemma 6.4. Let us remind that for a subspace $\mathcal{K}_{\Theta, \Delta, P}$, $\mathrm{Range} P(\xi) \perp \mathrm{Range} \Delta(\xi)$ a.e.

**Lemma 6.6.** *Assume that a system $\mathcal{K}_{\Theta_n, \Delta_n, P_n}$, $n \in \mathbb{N}$ be uniformly minimal. Let $\sigma_k$ denote a Borel support of $\Delta_k$, and $\tau_k$ be a Borel set where $\mathrm{Range} P_k(\xi) \oplus \mathrm{Range} \Delta_k(\xi) \neq E_*$. Then almost all points $\xi \in \mathbb{T}$ are covered by at most $d_* := \dim E_*$ sets $\sigma_k$ and by at most $d_*$ sets $\tau_k$. In particular, at most $d_*$ functions $\Theta_n \in H^\infty(E_n \to E)$ satisfy $\dim E_n > \dim E$.*

*Proof.* Let $\sigma_k$ denote a Borel support of $\Delta_k$. Consider all intersections of $d_* + 1$ sets $\sigma_k$ (there are countably many such intersections). If all such intersections have measure 0. the conclusion of the lemma is true. Suppose that there exists $d_* + 1$ sets $\sigma_k$ such that their intersection has positive measure. Without loss of generality we can assume that these are set $\sigma_1, \sigma_2, \ldots, \sigma_{d_*+1}$.

Take functions $f_k \in L^2(E_*)$ such that $f_k(\xi) \in \mathrm{Range} \Delta_k(\xi)$, $\|f_k(\xi)\| = 1$ on $\sigma_k$ and $f_k(\xi) = 0$ for $\xi \notin \sigma_k$.

Let us approximate functions $f_k$ uniformly by simple functions $\tilde{f}_k$, $\tilde{f}_k = \sum_j \chi_{\sigma_k^j} e_k^j$, where $\sigma_k^j$, $j = 1, 2, \ldots$ is a disjoint measurable covering of $\sigma_k$, $e_k^j \in E_*$, $\|e_k^j\| = 1$ and

$$\|f_k(\xi) - \tilde{f}_k(\xi)\| \leq \delta/(20(d_* + 1)) \qquad \text{a.e.,} \tag{6.1}$$

where $\delta$ is the constant of uniform minimality of the system $\mathcal{K}_{\Theta_n, \Delta_n, P_n}$, $n \in \mathbb{N}$.

Since for each $k = 1, 2, \ldots, d_* + 1$ sets $\sigma_k^j$ cover $\sigma_k$, one can find sets $\sigma_k^{j_k}$, $k = 1, 2, \ldots, d_* + 1$ with intersection $\sigma$ of positive Lebesgue measure. Define functions $g_k = m(\sigma)^{-1/2}\chi_\sigma e_k^{j_k} =: m(\sigma)^{-1/2}\chi_\sigma e_k$ (for the brevity of notation we skip the index $j_k$). Note that $\|g_k\| = 1$.

Lemma 6.5 implies that that for big enough $n \in \mathbb{N}$

$$\mathrm{dist}\left\{\overline{z}^n\begin{pmatrix} 0 \\ g_k \end{pmatrix}, \mathcal{K}_{\Theta_k,\Delta_k,P_k}\right\} \leq \delta/(10(d_*+1)), \qquad k = 1, 2, \ldots, d_* + 1,$$

where $\delta$ is the constant of uniform minimality of the system $\mathcal{K}_{\Theta_n,\Delta_n,P_n}$, $n \in \mathbb{N}$ (we can get as close to $\delta/(20(d_*+1))$ as we want, see (6.1)). We know that the system of subspaces $\mathcal{K}_{\Theta_k,\Delta_k,P_k}$ is uniformly minimal. Therefore, for big enough $n$ the system of functions $\overline{z}^n g_k = m(\sigma)^{-1/2}\overline{z}^n\chi_\sigma e_k$, $k = 1, 2, \ldots, d_*+1$ is uniformly minimal. That implies that the system of vectors $e_k \in E_*$, $k = 1, 2, \ldots, d_* + 1 = \dim E_* + 1$ is uniformly minimal, and that is impossible.

To prove the statement about $\tau_k$ let us first notice that if $f(\xi) \perp \mathrm{Range} P(\xi) \oplus \mathrm{Range} \Delta(\xi)$ a.e., then by Lemma 6.5

$$\begin{pmatrix} 0 \\ f \end{pmatrix} \in \mathcal{K}_{\Theta,\Delta,P}.$$

Suppose that there exist $d_* + 1 = \dim E_* + 1$ sets $\tau_k$ with intersection of $\tau$ positive measure. Again, as before, we can assume without loss of generality that these sets are $\tau_1, \tau_2 \ldots, \tau_{d_*+1}$. Take functions $f_k \in L^\infty(E_*)$ such that $f(\xi) \perp \mathrm{Range} P(\xi) \oplus \mathrm{Range}\Delta(\xi)$, $\|f(\xi)\| = 1$ a.e. on $\tau_k$, and $f_k(\xi) = 0$, $\xi \notin \tau_k$. Then we approximate functions $f_k$ uniformly by simple functions $\tilde{f}_k$, $\tilde{f}_k = \sum_j \chi_{\sigma_k^j} e_k^j$, where $\sigma_k^j$, $j = 1, 2, \ldots$ is a disjoint measurable covering of $\sigma_k$, $e_k^j \in E_*$, $\|e_k^j\| = 1$, and $\tilde{f}_k$ satisfy (6.1).

Again as above we find sets $\tau_k^{j_k}$, $k = 1, 2, \ldots, d_* + 1$, such that $\tau := \cap_{k=1}^{d_*+1}\tau_k^{j_k}$ has positive measure. Define $g_k := (m(\tau))^{-1/2}\chi_\tau e_k^{j_k} =: (m(\tau))^{-1/2}\chi_\tau e_k$. Then

$$\mathrm{dist}\left\{\begin{pmatrix} 0 \\ g_k \end{pmatrix}, \mathcal{K}_{\Theta_k,\Delta_k,P_k}\right\} \leq \delta/(10(d_*+1)), \qquad k = 1, 2, \ldots, d_* + 1,$$

and we can conclude that the system $g_k$, $k = 1, 2, \ldots, d_* + 1$ is uniformly minimal. So the vectors $e_k$ form a uniformly minimal system, what is impossible. $\square$

Let us remind that $P_{\Theta,\Delta,P}$ denote the orthogonal projection onto $\mathcal{K}_{\Theta,\Delta,P}$.

**Lemma 6.7.** *Let $g = \begin{pmatrix} 0 \\ f \end{pmatrix} \in \begin{pmatrix} H^2(E) \\ L^2(E_*) \end{pmatrix}$, and let the system of subspaces $\mathcal{K}_{\Theta_n,\Delta_n,P_n}$ be uniformly minimal. Then the following imbedding holds,*

$$\sum_n \|P_{\Theta_n,\Delta_n,P_n}g\|^2 \leq 2\dim E_* \cdot \|g\|^2, \qquad \forall g = \begin{pmatrix} 0 \\ f \end{pmatrix} \in \begin{pmatrix} H^2(E) \\ L^2(E_*) \end{pmatrix}.$$

*Proof.* Let $\sigma_k$, $\tau_k$ be the sets from Corollary 6.6. By Lemma 6.5

$$\|P_{\Theta_n,\Delta_n,P_n}g\|^2 = \|g\|^2 - (\|P_nf\|^2 + \|\mathbb{P}_+\Delta_n^*f\|^2) =$$
$$= \|f\|^2 - \|P_nf\|^2 - \|\mathbb{P}_+\Delta_n^*f\|^2 =$$
$$= \int_{\sigma_n\cup\tau_n} \|f(\xi)\|^2 dm(\xi) + \int_{\mathbb{T}\setminus(\sigma_n\cup\tau_n)} \|f(\xi)\|^2 dm(\xi) - \|P_nf\|^2 - \|\mathbb{P}_+\Delta_n^*f\|^2 \le$$
$$\le \int_{\sigma_n} \|f(\xi)\|^2 dm(\xi) + \int_{\tau_n} \|f(\xi)\|^2 dm(\xi) + \int_{\mathbb{T}\setminus(\sigma_n\cup\tau_n)} \|f(\xi)\|^2 dm(\xi)$$
$$- \int_{\mathbb{T}\setminus(\sigma_n\cup\tau_n)} \|P_n(\xi)f(\xi)\|^2 dm(\xi) = \int_{\sigma_n} \|f(\xi)\|^2 dm(\xi) + \int_{\tau_n} \|f(\xi)\|^2 dm(\xi)$$

Here in the last equality we used the fact that $\text{Range} P_n(\xi) = E_*$ on $\mathbb{T}\setminus(\sigma_n\cup\tau_n)$.

Taking the sum over all $n$ and applying Corollary 6.6 we get the required imbedding theorem. $\square$

## 7. Reduction to a scalar imbedding theorem

In this section we want to prove the following theorem. Recall that $P_{\Theta,\Delta,P}$ stands for the orthogonal projection onto $\mathcal{K}_{\Theta,\Delta,P}$

**Theorem 7.1.** *Let* $\dim E < \infty$, $\dim E_* < \infty$. *Suppose that the system of subspaces* $\mathcal{K}_{\Theta_n,\Delta_n,P_n}$ *is uniformly minimal and*

$$\sup_{\lambda\in\mathbb{D}}\sum_n(1-|\det\Theta_n(\lambda)|^2) =: C < \infty. \tag{7.1}$$

*Then the imbedding*

$$\sum_n \|P_{\Theta_n,\Delta_n,P_n}f\|^2 \le C'\|f\|^2 \tag{7.2}$$

*holds for all* $f \in \begin{pmatrix} H^2(E) \\ L^2(E_*) \end{pmatrix}$ *and, moreover,*
$C' = C'(C,\dim E,\dim E_*)$.

First of all let us remind that by Lemmas 6.4 and 6.6 at most $\dim E + \dim E_*$ matrices $\Theta_n$ are non-square, so one always can omit this non-square terms without influencing the imbedding theorem.

If we would be able to show that the uniform minimality of the system of subspaces $\mathcal{K}_{\Theta_n,\Delta_n,P_n}$ implies the condition (7.1) the main result would be proved. Indeed, Theorem 7.1 then implies the imbedding (7.2), which is exactly the first imbedding in Theorem 3.1. (Formally, for the imbedding 1 in Theorem 3.1 we need only (7.2) for $f \in \text{span}\{\mathcal{K}_{\Theta_n,\Delta_n,P_n} : n = 1,2,\ldots\}$, but of course it is the same.)

The dual system to $\mathcal{K}_{\Theta_n,\Delta_n,P_n}$, $n = 1, 2, \ldots$ is a system of invariant subspaces of $T^*$, and $T^*$ is a completely non-unitary contraction with finite defects. It can be represented as (unitarily equivalent to) a system of subspaces $\mathcal{K}_{\Theta'_n,\Delta'_n,P'_n}$ in a space $\begin{pmatrix} H^2(E_*) \\ L^2(E) \end{pmatrix}$, where $\Theta'_n, \Delta'_n, P'_n$ are some operator valued functions. It is possible to write some formulas to relate $\Theta_n, \Delta_n, P_n$ and $\Theta'_n, \Delta'_n, P'_n$, but we do not need that!. We know that the dual system is of course also uniformly minimal, so we can apply Theorem 7.1 without even knowing anything about functions $\Theta'_n, \Delta'_n, P'_n$. We get another imbedding from Theorem 3.1 and we are done.

To prove Theorem 7.1 we need the following result. Let $\mathbb{P}_-$ denotes the orthogonal projection from $L^2$ onto $H^2$.

**Theorem 7.2.** *Let $F_n$ be a sequence of functions, $F_n \in H^\infty$. Then the imbedding*
$$\sum_n \|\mathbb{P}_-\overline{F}_n f\|_2^2 \leq C\|f\|^2 \qquad \forall f \in H^2$$

*holds if and only if*
$$\sup_{\lambda \in \mathbb{D}} \sum_n (|F|^2(\lambda) - |F(\lambda)|^2) =: C' < \infty; \qquad (7.3)$$

*here $|F|^2(\lambda)$ denote the harmonic extension of $|F|\big|\mathbb{T}$ at the point $\lambda \in D$. Moreover, the constants $C$ and $C'$ are equivalent, i.e. $A^{-1}C \leq C' \leq AC$ where $A$ is some absolute constant.*

*Proof.* Consider a Hankel operator $H_{\overline{F}} : H^2 \to H^2_-(\ell^2) = H^2_- \otimes \ell^2$
$$(H_{\overline{F}} f)_k := H_{\overline{F}_k} f = \mathbb{P}_- F_k f.$$

In other words the symbol $\overline{F}$ of the Hankel operator $H_{\overline{F}}$ is the column with entries $\overline{F}_k$, $k = 1, 2, \ldots$

Then the sum in (7.3) is the square of an equivalent BMO norm (the so-called Garsia norm) of the function $\mathbb{P}_-\overline{F}$. It is well known that for Hankel operators whose symbol is a column $\overline{F}$ the norm of the corresponding Hankel operator and the BMO norm of $\mathbb{P}_-\overline{F}$ are equivalent. For scalar case this fact is contained in all textbooks, and for general case the proof follows the lines there. Note that the result does not hold for Hankel operators with general (operator-valued) symbols. It is essential that our symbol has only one column. □

**Corollary 7.3.** *Let $F_n$ be a sequence of operator valued functions, with*
$$F_n \in H^\infty(E_n \to E).$$

*Then the imbedding*
$$\sum_n \|\mathbb{P}_- F_n^* f\|_2^2 =: C\|f\|_2^2 \qquad \forall f \in H^2(E)$$

*holds if and only if*

$$\sup_{e \in E, \|e\|=1} \sup_{\lambda \in \mathbb{D}} \sum_n (\|F_n^* e\|^2(\lambda) - \|F_n(\lambda)^* e\|^2) =: C' < \infty; \tag{7.4}$$

*here $\|F_n^* e\|^2(\lambda)$ denote the harmonic extension of the function $\|F_n^* e\|^2$ defined the unit circle $\mathbb{T}$ to a point $\lambda \in \mathbb{D}$. Moreover, the constants $C$ and $C'$ are equivalent, i.e. $A^{-1}C \leq C' \leq AC$ where $A$ depends only on $\dim E$.*

*Proof.* Consider a restriction of the imbedding to the functions $f$ of form $f = \varphi e$ where $\varphi$ is a scalar $H^2$ function and $e \in E$ is a fixed vector, $\|e\| = 1$. The condition (7.3) guarantee that the imbedding holds on that functions. Consider an orthonormal basis $\{e_n\}_{n=1}^d$, $d := \dim E$ in $E$. Then the space $H^2(E)$ is an orthogonal sum of subspaces $\mathcal{H}_k$, $\mathcal{H}_k$ consists of all functions of form $f = \varphi e$ where $\varphi$ is a scalar $H^2$ function. We have imbedding on each $\mathcal{H}_k$ and since there are only finitely many of them, the imbedding holds on $H^2(E)$. □

*Proof of Theorem 7.1.* By Lemma 6.7 it is enough to check the imbedding on vectors $\begin{pmatrix} f \\ 0 \end{pmatrix} \in \begin{pmatrix} H^2(E) \\ L^2(E) \end{pmatrix}$ By Corollary 6.2

$$\left\| P_{\Theta, \Delta, P} \begin{pmatrix} f \\ 0 \end{pmatrix} \right\|^2 = \|f\|^2 - \|\mathbb{P}_+ \Theta^* f\|^2 \leq$$
$$\leq \|(I - \Theta\Theta^*)^{1/2} f\|^2 + \|\Theta^* f\|^2 - \|\mathbb{P}_+ \Theta^* f\|^2 =$$
$$= \|(I - \Theta\Theta^*)^{1/2} f\|^2 + \|\mathbb{P}_- \Theta^* f\|^2$$

So

$$\sum_n \left\| P_{\Theta_n, \Delta_n, P_n} \begin{pmatrix} f \\ 0 \end{pmatrix} \right\|^2 \leq \sum_n \|(I - \Theta_n \Theta_n^*)^{1/2} f\|^2 + \sum_n \|\mathbb{P}_- \Theta_n^* f\|^2 \tag{7.5}$$

First of all, by Lemmas 6.4 and 6.6 at most $\dim E + \dim E_*$ matrices $\Theta_n$ are non-square. So we can assume that in the sum we have only square matrices, in other words that $\Theta_n \in H^\infty(E \to E)$ and $\det \Theta_n$ are well defined.

The condition (7.1) clearly implies (because $\Theta_n(\xi)$, $\xi \in \mathbb{T}$ are contractions)

$$\sup_{e \in E, \|e\|=1} \sup_{\lambda \in \mathbb{D}} \sum_n (1 - \|\Theta_n(\lambda)^* e\|^2) < \infty,$$

and the last condition implies

$$\sup_{e \in E, \|e\|=1} \sup_{\lambda \in \mathbb{D}} \sum_n (\|\Theta_n^* e\|^2(\lambda) - \|\Theta_n(\lambda)^* e\|^2) < \infty.$$

By Corollary 7.3 that implies

$$\sum_n \|\mathbb{P}_- \Theta_n^* f\|^2 \leq C \|f\|^2 \qquad \forall f \in H^2(E),$$

so we estimated the second sum in (7.5).

To estimate the first sum, let us denote by $\sigma_n$ a Borel support of $\Delta_n$, see Corollary 6.6. For $\xi \in \mathbb{T} \setminus \sigma_n$ operators $\Theta_n(\xi)$ are isometries, and since we assumed $\Theta_n \in H^\infty(E \to E)$, $\dim E < \infty$, it follows that $(I - \Theta_n(\xi)\Theta_n(\xi)^*) = 0$ for $\xi \in \mathbb{T} \setminus \sigma_n$. Therefore

$$\|(I - \Theta_n \Theta_n^*)^{1/2} f\|^2 \leq \int_{\sigma_n} \|f(\xi)\|^2 dm(\xi).$$

By Corollary 6.6 almost all points of $\mathbb{T}$ are covered by at most $\dim E_*$ sets $\sigma_n$ and therefore

$$\sum_n \|(I - \Theta_n \Theta_n^*)^{1/2} f\|^2 \leq \dim E_* \|f\|^2$$

so we estimated the first sum in (7.5). $\square$

So we made one more reduction. Now to prove the main result (Theorem 2.1) it remains to show that the uniform minimality of the system of subspaces $\mathcal{K}_{\Theta_n, \Delta_n, P_n}$, $n = 1, 2, \ldots$ implies the condition (7.1)

## 8. Geometry of vector reproducing kernels.

Consider in a Hilbert space $\mathcal{H}$ a system of subspaces $\mathcal{E}_n$, $n = 1, 2, \ldots$. Suppose that the system is a Riesz basis, and let $R$ be a normalized orthogonalizer for this system, see Section 1.1. Let us remind that an orthogonalizer is a bounded invertible operator that maps the system $\mathcal{E}_n$, $n = 1, 2, \ldots$ into an orthogonal one, and the "normalized" mean that $R|\mathcal{E}_n$ is an isometry for all $n$.

It follows from the definition of the normalized orthogonalizer that orthogonalizer that for any (finite) set of vectors $f_n \in \mathcal{E}_n$

$$\|R\|^{-2} \sum_n \|f_n\|^2 \leq \Big\|\sum_n f_n\Big\|^2 \leq \|R^{-1}\|^2 \sum_n \|f_n\|^2. \tag{8.1}$$

We need the following simple lemma.

**Lemma 8.1.** *Let a system of subspaces $\mathcal{E}_n$, $n = 1, 2, \ldots$ in a Hilbert space $\mathcal{H}$ be a Riesz basis, and let $R$ be its normalized orthogonalizer. Then for a Hilbert space $E$ the system $\mathcal{E}_n \otimes E$, $n = 1, 2, \ldots$ in $\mathcal{H} \otimes E$ is a Riesz basis and $R \otimes I_E$ is a normalized orthogonalizer for this system.*

We apply this Lemma in the following situation. Consider a countable subset $\sigma \subset \mathbb{D}$, and let $\{k_\lambda\}_{\lambda \in \sigma}$ be a system of normalized reproducing kernels in $H^2$,

$$k_\lambda(z) := \frac{(1 - |\lambda|^2)^{1/2}}{1 - \overline{\lambda} z},$$

$\|k_\lambda\|_2 = 1$.

Suppose the system $\{k_\lambda\}_{\lambda \in \sigma}$ is a Riesz basis (i.e. the system of corresponding one-dimensional subspaces is a Riesz basis), and let $R$ be its (normalized) orthogonalizer.

In our case it is simply an operator that maps the system $\{k_\lambda\}_{\lambda\in\sigma}$ into orthonormal system. The following result is a simple corollary of Lemma 8.1

**Lemma 8.2.** *Let the system $\{k_\lambda\}_{\lambda\in\sigma}$ be a Riesz basis, let $R$ be its orthogonalizer, and Let $E$ be a Hilbert space. Then for any finite set of vectors $f_\lambda \in E$*

$$\|R\|^{-2} \sum_{\lambda\in\sigma} \|f_\lambda\|^2 \leq \Big\| \sum_{\lambda\in\sigma} k_\lambda f_\lambda \Big\|^2 \leq \|R^{-1}\|^2 \sum_{\lambda\in\sigma} \|f_\lambda\|^2.$$

*Proof.* The statement of the lemma is simply the inequality (8.1), because by Lemma 8.1 $R\otimes I$ is a normalized orthogonalizer of the system $\{\operatorname{span}\{k_\lambda\}\otimes E\}_{\lambda\in\sigma}$, and clearly $\|R\otimes I\| = \|R\|$, $\|(R\otimes I)^{-1}\| = \|R^{-1}\|$. □

We need several well known facts about geometry of a system of (scalar) reproducing kernels. The following theorem can be found for example in [1].

**Theorem 8.3.** *The system of reproducing kernels $\{k_\lambda\}_{\lambda\in\sigma}$ is a Riesz basis if and only if the set $\sigma$ is an interpolating set,*

$$\inf_{\lambda\in\sigma} \prod_{\mu\in\sigma\setminus\{\lambda\}} |b_\mu(\lambda)| =: \delta > 0$$

*where $b_\lambda$ is a Blaschke factor, $b_\lambda(z) \stackrel{\text{def}}{=} (|\lambda|/\lambda)(\lambda - z)(1 - \overline{\lambda}z)^{-1}$.*

*Moreover, the measure of nonorthogonality $\|R\| \cdot \|R^{-1}\|$ of the system ($R$ is a normalized orthogonalizer) admits the estimate above by a constant $C(\delta)$ depending only on $\delta$.*

Let $\theta_n$ be scalar valued inner functions in $H^\infty$, and let $\mathcal{K}_{\theta_n}$ denote the corresponding $S^*$-invariant subspace, $\mathcal{K}_{\theta_n} := H^2 \ominus \theta_n H^2$.

We need the following theorem that was proved by Vasyunin in [10], see also [1, Lecture IX]

**Theorem 8.4** (Vasyunin). *Let the system of subspaces $\mathcal{K}_{\theta_n}$, where $\theta_n$ are scalar inner functions in $H^\infty$, be uniformly minimal and let $\delta$ be its constant of uniform minimality. Then the system is a Riesz basis, and there exists a constant $CV(\delta)$ such that the measure of non-orthogonality $\|R\| \cdot \|R^{-1}\|$ (recall that $R$ is the normalized orthogonalizer of the system) admits the following estimate*

$$\|R\| \cdot \|R^{-1}\| \leq CV(\delta)$$

*Remark 8.5.* The above theorem can be also proved using Theorem 5.1 above. To do this it is enough to notice that the uniform minimality of the system of subspaces $\mathcal{K}_{\theta_n}$ implies

$$\inf_n \inf_{\lambda\in\mathbb{D}} \left\{ |\theta_n(\lambda)| + \prod_{k\neq n} |\theta_k(\lambda)| \right\} =: \delta' > 0$$

and that $\delta' \leq 2\delta$, see [1, Lecture IX]. The above condition implies (5.2) from Theorem 5.1, and therefore the imbedding (5.1). The dual system has the same geometry, so the imbedding for the dual system holds automatically.

8.1. **Carleson measures and bases of reproducing kernels.** For an arc $I \subset \mathbb{T}$ let $S(I)$ denote the *Carleson square*,
$$S(I) := \{z \in \mathbb{D} \,:\, z/|z| \in I,\, 1 - |I| \leq |z| < 1\}.$$
A measure $\mu$ on the disk $\mathbb{D}$ is called Carleson if
$$\sup_{I \subset \mathbb{T}} \frac{\mu(S(I))}{|I|} =: K < \infty;$$
here supremum is taken over all arcs $I \subset \mathbb{T}$.

The constant $K$ is called the *Carleson norm* of the measure $\mu$.

The following theorem is well known, see for example monographs [15, 1].

**Theorem 8.6.** *Let $\mu$ be a measure in $\mathbb{D}$. The following statements are equivalent*

1. *The measure $\mu$ is Carleson;*
2. $\sup_{\lambda \in \mathbb{D}} \int |k_\lambda(z)|^2 d\mu(z) =: K' < \infty$, *where $k_\lambda$ is the normalized reproducing kernel of $H^2$, $k_\lambda(z) = (1 - |\lambda|^2)^{1/2}(1 - \overline{\lambda}z)^{-1}$;*
3. *the imbedding*
$$\int_{\mathbb{D}} |f(z)|^2 d\mu(z) \leq K'' \|f\|^2, \qquad f \in H^2$$
    *holds.*

*Moreover, the Carleson norm $K$ of $\mu$, the constant $K'$ and the best possible $K''$ are equivalent in the sense of two-sided estimates with some absolute constants.*

The following description of interpolating sets is well known an again can be found for example in [15, 1]

**Theorem 8.7.** *The set $\sigma \subset \mathbb{D}$ satisfies the Carleson interpolation condition*
$$\inf_{\lambda \in \sigma} \prod_{\mu \in \sigma \setminus \{\lambda\}} |b_\mu(\lambda)| =: \delta > 0$$
*if and only if*

1. *The points of $\sigma$ are separated in hyperbolic metric, i.e.*
   $\inf \{|b_\lambda(\mu)| \,:\, \lambda, \mu \in \sigma,\, \lambda \neq \mu\} =: \alpha > 0$, *and*
2. *the measure $\sum_{\lambda \in \sigma}(1 - |\lambda|^2)\delta_\lambda$, where $\delta_\lambda$ stands for the unit mass at $\lambda$, is Carleson.*

*Moreover the Carleson constant $\delta$ can be estimated by the constant depending only on $\alpha$ and the Carleson norm $K$ of the measure*
$$\delta \geq \delta(\alpha, K),$$

*and vice versa $K$ and $\rho$ can be estimated by constant depending only on $\delta$,*

$$K \leq K(\delta), \qquad \alpha \geq \alpha(\delta)$$

In [15] the estimates are not stated, but they are contained in the proofs.

## 9. Carleson contours

We need the following result

**Theorem 9.1.** *Given $\varepsilon > 0$ there exists $\varepsilon' = \varepsilon'(\varepsilon) > 0$, $\varepsilon' < \varepsilon$ such that for any $\varphi \in H^\infty$, $\|f\|_\infty \leq 1$ one can find a region $\mathcal{O}$ satisfying*

$$\{z \in \mathbb{D} : |\varphi(z)| < \varepsilon'\} \subset \mathcal{O} \subset \{z \in \mathbb{D} : |\varphi(z)| \leq \varepsilon\}$$

*and such that the arc-length on $\gamma := \partial \mathcal{O} \cap \mathbb{D}$ is a Carleson measure with the Carleson norm bounded by an absolute constant $K$.*

This theorem was proved in [16]. More precisely it was proved only for inner functions, but the same proof works for general $H^\infty$ functions as well. There was also a slight inaccuracy in the proof, so for the sake of completeness, and to make the papers self-contained, we present the proof in the Appendix.

The boundary $\partial \mathcal{O} \cap \mathbb{D}$ of $\mathcal{O}$ is usually called a Carleson contour for $\varphi$. There exists quite a few constructions of such contours. We chose the Bourgain's construction because for our problem it is essential the the Carleson norm of the arc-length on $\partial \mathcal{O}$ does not depend on $\varepsilon$ (at least this makes our life a lot easier).

## 10. The main construction

Let us pick a small $\varepsilon$ and for each $n$ construct a Carleson contour for $\det \Theta_n$. The constant $\varepsilon$ will be chosen later. For now we only assume that it is less or equal than $\varepsilon$ from Lemma 6.4. As it was said above we can assume without loss of generality that all matrix-functions $\Theta_n$ are $d \times d$ matrix-functions, so determinants are well defined.

So we have open sets $\mathcal{O}_n$ with rectifiable boundary $\partial \mathcal{O}_n$, such that

$$\{z \in \mathbb{D} : |\det \Theta_n(z)| < \varepsilon'^d\} \subset \mathcal{O}_n \subset \{z \in \mathbb{D} : |\Theta_n(z)| \leq (\varepsilon)^d\}$$

and the arc-length on $\gamma_n := \partial \mathcal{O}_n \cap \mathbb{D}$ is a Carleson measure with Carleson norm at most $K$.

Let us recall that the hyperbolic distance between two points $\lambda$ and $\mu$ in $\mathbb{D}$ is defined as $\rho(\lambda, \mu) := \frac{1}{2} \log\{(1 + |b_\lambda(\mu)|)/(1 - |b_\lambda(\mu)|)\}$, where $b_\lambda$ denotes a Blaschke factor with zero at $\lambda$. It is well known, see [11, 12] that $\rho$ is indeed a metric, and it is easy to see that if $|b_\lambda(\mu)|$ is small, then $\rho(\lambda, \mu) \approx |b_\lambda(\mu)|$.

Fix a small constant $\alpha$, $0 < \alpha < 0.1$. For each $n$ pick on the contour $\gamma_n$ a discrete set $\sigma_n$ such that the points in this set are hyperbolically separated, $|b_\lambda(\mu)| > \alpha > 0$, $\lambda, \mu \in \sigma_n$, $\lambda \neq \mu$, but are dense enough, namely that for any point $z \in \gamma_n$ one can find a point $\lambda \in \sigma_n$ for which $|b_\lambda(z)| < \alpha$.

It can be done by putting points $\lambda$ on the contour in such a fashion that hyperbolic discs of (hyperbolic) radius $\rho = \frac{1}{2}\log\{(1+\alpha)/(1-\alpha)\}$ centered at that points cover the contour $\gamma_n$, but that every next point (center) is chosen outside the union of hyperbolic disk of radius $\rho$ with centers in the points already constructed.

To show that we indeed can cover the contour this way, let us split the disk $\mathbb{D}$ into "layers" $D_m$, $D_m := \{z \in \mathbb{D} : 1 - 2^{-m} \leq |z| < 1 - 2^{-m-1}\}$, $m = 0, 1, 2, \ldots$, and putting points $\lambda$ on the contour first in the "layer" $D_0$ while it is possible, then in the "layer" $D_2$, and so on. By the construction the hyperbolic disk of radius $\rho/2$ centered at points $\lambda \in \sigma_n$ are disjoint. Since each "layer" $D_m$ has a finite hyperbolic area, and the hyperbolic area of a hyperbolic disk depends only of radius, we can put only finitely many points of $\sigma_n$ in each "layer". So it is indeed possible to cover all the contour $\gamma_n$ by countably many hyperbolic disks.

Note that the measure $\sum_{\lambda \in \sigma_n}(1 - |\lambda|^2)\delta_\lambda$ is Carleson, and its Carleson norm can be estimated by the Carleson norm of $|dz|\big|\gamma_n$ (arc-length on $\gamma_n$),

$$A(\alpha) \left\| |dz|\big|\gamma_n \right\|_{\mathrm{Carl}}^2 \leq \left\| \sum_{\lambda \in \sigma_n}(1 - |\lambda|^2)\delta_\lambda \right\|_{\mathrm{Carl}} \leq B(\alpha) \left\| |dz|\big|\gamma_n \right\|_{\mathrm{Carl}}$$

By the construction points of $\sigma_n$ are hyperbolically separated, so by Theorem 8.7 the set $\sigma_n$ is an interpolating set, and its interpolation constant $\delta = \delta(\sigma_n)$ can be estimated by a constant depending only on $\alpha$,

$$\inf_{\lambda \in \sigma_n} \prod_{\mu \in \sigma_n \setminus \{\lambda\}} |b_\mu(\lambda)| =: \delta(\sigma_n) \geq \delta(\alpha) > 0$$

Therefore by Theorem 8.3 the system $\{k_\lambda\}_{\lambda \in \sigma_n}$ is a Riesz basis, and its measure of non-orthogonality $\|R\| \cdot \|R^{-1}\|$ admits the estimate

$$\|R\| \cdot \|R^{-1}\| \leq C = C(\alpha). \tag{10.1}$$

Let $B_n$ be a Blaschke product with zero set $\sigma_n$,

$$B_n := \prod_{\lambda \in \sigma_n} b_\lambda.$$

Note that for $B_n$ the following estimate holds

$$|B_n(z)| < \alpha, \qquad z \in \gamma_n,$$

because by the construction $|b_\lambda(z)| < \alpha$ for some $\lambda \in \sigma_n$.

The following lemma allows us to replace $\det \Theta_n$ by the constructed above Blaschke products $B_n$.

**Lemma 10.1.** *If*

$$\sup_{\lambda \in \mathbb{D}} \sum_n (1 - |B_n(\lambda)|^2) =: K < \infty \tag{10.2}$$

*then the condition*

$$\sup_{\lambda \in \mathbb{D}} \sum_n (1 - |\det \Theta_n(\lambda)|^2) \leq C = C(d, d_*, K, \alpha, \varepsilon, \varepsilon') < \infty. \tag{10.3}$$

*holds.*

*Proof.* Let $h_n$ be outer functions in $H^\infty$ such that

$$|h_n(\xi)| = \max\{|\det \Theta_n(\xi)|, \varepsilon'^d\}, \qquad \xi \in \mathbb{T}.$$

Clearly $\|h_n\|_\infty \leq 1$.

Notice that for big enough $N = N(d, \alpha, \varepsilon')$

$$|B_n^N(z)| < (\varepsilon')^d \leq |\det \Theta_n(z)|, \qquad z \in \gamma_n.$$

On the other hand on for the points $\xi \in \mathbb{T}$ that are angular limits of points in $\mathbb{D} \setminus \mathcal{O}_n$ the inequality $|\det \Theta_n(\xi)| \geq \varepsilon'^d$ holds, thus $|\det \Theta_n(\xi)| \geq |h_n(\xi)|$ for almost all such points. So by maximum modulus principle (we have to use a generalization of a classical maximum modulus principle, see Lemma 10.2 below)

$$|h_n(z)| \cdot |B_n(z)|^N \leq |\det \Theta_n(z)|, \qquad z \in \mathbb{D} \setminus \mathcal{O}_n.$$

Therefore for a point $z \in \mathbb{D}$

$$\sum_{n : z \notin \mathcal{O}_n} (1 - |\det \Theta_n(z)|^2) \leq \sum_n (1 - |h_n(z)|^2 |B_n(z)|^{2N}) \leq$$

$$\leq \sum_n (1 - |h_n(z)|^2) + \sum_n (1 - |B_n(z)|^{2N}) \leq$$

$$\leq \sum_n (1 - |h_n(z)|^2) + N \sum_n (1 - |B_n(z)|^2).$$

Here the second inequality follows from a trivial estimate

$$1 - \alpha_1 \alpha_2 = \alpha_2(1 - \alpha_1) + (1 - \alpha_2) \leq (1 - \alpha_1) + (1 - \alpha_2), \qquad 0 \leq \alpha_1, \alpha_2 \leq 1$$

By Lemma 6.4 any point $z \in \mathbb{D}$ is covered by at most $d = \dim E$ sets $\{z \in \mathbb{D} : |\det \Theta_n(z)| < \varepsilon^d\} \supset \mathcal{O}_n$ (we assume that $\varepsilon$ is small enough), hence

$$\sum_n (1 - |\det \Theta_n(z)|^2) \leq \sum_n (1 - |h_n(z)|^2) + N \sum_n (1 - |B_n(z)|^2) + d.$$

To complete the proof of the lemma we need to estimate $\sum_n (1 - |h_n(z)|^2)$. Let $s_k$ denote a Borel subset of $\mathbb{T}$ such that $|h_k(\xi)| < 1$, $\xi \in s_k$ and $|h_k(\xi)| = 1$, $\xi \in \mathbb{T} \setminus s_k$. By Corollary 6.6 almost all points $\xi \in \mathbb{T}$ are covered by at most $d_* = \dim E_*$ sets $s_k$. Therefore the product $\prod_n |h_n|$ converges on $\mathbb{T}$ and moreover

$$\prod_n |h_n(\xi)| \geq (\varepsilon')^{d_*}, \qquad \xi \in \mathbb{T}.$$

Therefore, if we normalize all $h_n$ by the condition $h_n(0) > 0$, then clearly the product $\prod_n h_n$ converges in $\mathbb{D}$ and moreover, by maximum modulus principle

$$\left|\prod_n h_n(z)\right| \geq (\varepsilon')^{d_*}, \qquad z \in \mathbb{D}.$$

It follows that

$$(\varepsilon')^{d_*} \leq \prod_n |h_n(z)| = \exp\left(\frac{1}{2}\sum_n \log |h_n(z)|^2\right) \leq$$

$$\leq \exp\left(-\frac{1}{2}\sum_n (1 - |h_n(z)|^2)\right).$$

Therefore

$$\sum_n (1 - |h_n(z)|^2) \leq 2d_* \log \frac{1}{\varepsilon'}.$$

□

The following lemma is the maximum modulus principle we have used in the proof of Lemma 10.1.

**Lemma 10.2.** *Let $\Omega \subset \mathbb{D}$ be an open set, $f$ be a bounded analytic in $\Omega$ function, continuous on $\operatorname{clos}\Omega \cap \mathbb{D}$. Suppose that $|f| \leq 1$ on $\mathbb{D} \cap \partial\Omega$ and $\limsup_{z\to\xi}|f(z)| \leq 1$ for all points $\xi \in \mathbb{T}$ that are angular limits of points in $\Omega$ (lim sup also means angular upper limit). Then $|f| \leq 1$ on $\Omega$.*

Such results are for sure familiar to specialists. The lemma stated above (even stringer version) can be found for example in [14], see Lemma 63 there.

10.1. **Choice of $\varepsilon$.** Let $\delta$ be the constant of uniform minimality of the system $\mathcal{E}_n = \mathcal{K}_{\Theta_n, \Delta_n, P_n}$, $n = 1, 2, \ldots$, and let $CV(\delta/2)$ be the constant from Theorem 8.4. Let $C = C(\alpha)$ be the constant from (10.1). Recall that we fixed the small constant $\alpha$. We will chose $\varepsilon$ smaller than $\varepsilon$ from Lemma 6.4 and such that

$$\varepsilon C \cdot CV(\delta/2) < \delta/10 \tag{10.4}$$

The choice will be clear from what follows below.

10.2. **Reduction to "almost scalar" case.** Consider a point $\lambda \in \sigma_n$. Since by the construction $|\det \Theta_n(\lambda)| < \varepsilon^d$, there exists a vector $e_\lambda \in E$, $\|e_\lambda\| = 1$, such that $\|\Theta(\lambda)^* e_\lambda\| < \varepsilon$.

Consider a system of vector reproducing kernels $\{k_\lambda e_\lambda\}_{\lambda \in \sigma}$ in $H^2(E)$, where $k_\lambda(z) = (1 - |\lambda|^2)^{1/2}(1 - \overline{\lambda}z)^{-1}$ is the normalized reproducing kernel in $H^2$.

Let $\varepsilon$ be the constant chosen above in Section 10.1. Let $e^1, e^2, \ldots, e^N$ be a finite $\varepsilon$-net for the unit sphere in $E$.

For each $n$ split the zero set $\sigma_n$ of $B_n$ into $N$ disjoint sets $\sigma_n^1, \sigma_n^2, \ldots, \sigma_n^N$ such that for any $\lambda \in \sigma_k$ the corresponding vector $e_\lambda$ is close to $e^k$

$$\|e_\lambda - e^k\| < \varepsilon, \qquad \lambda \in \sigma_n^k. \tag{10.5}$$

Let $B_n^k$ be a Blaschke product with zero set $\sigma_n^k$ (if $\sigma_n^k = \varnothing$ we put $B_n^k := 1$). Clearly $B_n = \prod_{k=1}^N B_n^k$.

**Lemma 10.3.** *For Blaschke products $B_n$ and $B_n^k$ constructed above*

$$\sup_{\lambda \in \mathbb{D}} \sum_n (1 - |B_n(\lambda)|^2) \leq \sum_{k=1}^N \sup_{\lambda \in \mathbb{D}} \sum_n (1 - |B_n^k(\lambda)|^2).$$

*Proof.* This lemma admits a simple geometrical interpretation and a simple proof in terms of imbedding theorems, but we will use the following elementary inequality:

$$1 - \alpha_1 \alpha_2 = \alpha_2(1 - \alpha_1) + (1 - \alpha_2) \leq (1 - \alpha_1) + (1 - \alpha_2),$$

for $0 < \alpha_1, \alpha_2 \leq 1$. By induction

$$1 - \alpha_1 \alpha_2 \ldots \alpha_n \leq (1 - \alpha_1) + (1 - \alpha_2) + \ldots + (1 - \alpha_n),$$

for $0 < \alpha_1, \alpha_2, \ldots, \alpha_n \leq 1$. $\square$

## 11. The main Lemma

**Lemma 11.1.** *Let $\mathcal{K}_{\Theta_n, \Delta_n, P_n}$, $n = 1, 2, \ldots$ be a uniformly minimal system of subspaces in $\begin{pmatrix} H^2(E) \\ L^2(E_*) \end{pmatrix}$, and let $\delta$ be its constant on uniform minimality. Let $B_n$ be interpolating Blaschke products with zero sets $\sigma_n$ and let for all $n$ the measure of non-orthogonality (see Introduction) $\|R\| \cdot \|R^{-1}\|$ of the family of reproducing kernels $\{k_\lambda\}_{\lambda \in \sigma_n}$ is at most $C$. Let $\varepsilon > 0$ be a small such that*

$$C \cdot CV(\delta/2) \cdot \varepsilon < \delta/10,$$

*where $CV$ is the constant from Theorem 8.4. Suppose that there exists a vector $e \in E$, $\|e\| = 1$ such that for all $n$ and for all $\lambda \in \sigma_n$*

$$\|\Theta_n(\lambda)^* e\| < 2\varepsilon.$$

*Then the system of subspaces $\mathcal{K}_{B_n} \subset H^2$ (scalar), $n = 1, 2, \ldots$ is uniformly minimal with constant of uniform minimality at least $\delta/2$.*

This lemma clearly implies the main result. For a fixed $k$ the constructed above Blaschke products $B_n^k$ and matrix functions $\Theta_n$ satisfy the assumptions of the lemma with $\sigma_n := \sigma_n^k$, $e := e^k$ because

$$\|\Theta_n(\lambda)^* e^k\| \leq \|\Theta_n(\lambda)^* e_\lambda\| + \|\Theta_n(\lambda)^*(e_\lambda - e^k)\| < 2\varepsilon$$

for $\lambda \in \sigma_n^k$. Therefore the system of subspaces $\mathcal{K}_{B_n^k}$ satisfy the Carleson–Vasyunin condition (CV). This implies

$$\sum_n (1 - |B_n^k(\lambda)|^2) \leq K = K(\delta) < \infty, \qquad \lambda \in \mathbb{D}.$$

As it was shown above in Section 10 it is enough to prove the main theorem.

11.1. **Idea of the proof of the main lemma.** The idea is very simple and looks very naïve. Suppose that the constant of the uniform minimality of the system $\mathcal{K}_{B_n}$, $n = 1, 2, \ldots$ is less than $\delta/2$. That means that for a vector $f_n \in \mathcal{K}_{B_n}$, $\|f_n\| = 1$ there exist vectors $f_k \in \mathcal{K}_{B_k}$ such that $\|f_n - \sum_{k \neq n} f_k\| < \delta/2$.

Consider the spaces $\mathcal{K}_{B_n} \otimes e \subset H^2(E)$ imbedded in $\begin{pmatrix} H^2(E) \\ L^2(E_*) \end{pmatrix}$. It can be shown that the subspaces $\mathcal{K}_{B_n}$ are close to $\mathcal{K}_{\Theta_n, \Delta_n, P_n}$, and we are going to obtain a similar approximation for some $\tilde{f}_n \in \mathcal{K}_{\Theta_n, \Delta_n, P_n}$, $\|\tilde{f}_n\| = 1$

$$\left\| \tilde{f}_n - \sum_{k \neq n} \tilde{f}_k \right\| < \delta, \qquad f_k \in \mathcal{K}_{\Theta_k, \Delta_k, P_k},$$

that contradicts the fact that the constant of uniform minimality of the original system is $\delta$.

The main trick here is that one has to take not too many summands $\tilde{f}_k$ in the approximation, so that the accumulated "round-off" error caused by the deviation from the scalar case does not worsen the estimate obtained for this case by more that $\delta/2$. Quantitatively this "not too many" is expressed by the following lemma, see [8, 7, 9].

Here and in what follows $\delta(\{\mathcal{E}_n\})$ denotes the constant of uniform minimality of the system of subspaces $\{\mathcal{E}_n\}$. We assume that $\delta(\{\mathcal{E}_n\}) = 0$ if the system is not uniformly minimal.

**Lemma 11.2.** *Let $\delta(\{\mathcal{E}_n\}) < \delta < 1$. There exists a finite subset $\mathcal{N}$ of indices such that*

$$\delta(\{\mathcal{E}_n\}_{n \in \mathcal{N}}) < \delta \tag{11.1}$$

*but for any $k \in \mathcal{N}$*

$$\delta(\{\mathcal{E}_n\}_{n \in \mathcal{N} \setminus \{k\}}) \geq \delta. \tag{11.2}$$

*Proof.* Chose a finite set of indices $\mathcal{N}$ satisfying (11.1) (such choice is always possible because $\delta(\{\mathcal{E}_n\}) < \delta$. If (11.2) holds, the lemma is proved. Otherwise remove from $\mathcal{N}$ a point $k$ for which $\delta(\{\mathcal{E}_n\}_{n \in \mathcal{N} \setminus \{k\}}) < \delta$, check (11.2) again, and so on. Obviously, we eventually obtain the desired set $\mathcal{N}$, since the constant of uniform minimality of a system consisting of one subspace is 1. $\square$

*Proof of the main lemma (Lemma 11.1).* Suppose the conclusion of the lemma does not hold and $\delta(\{\mathcal{K}_{B_n}\}_{n\in\mathbb{N}}) < \delta/2$. Let us apply Lemma 11.2 with constant $\delta/2$ to the system $\{\mathcal{K}_{B_n}\}_{n\in\mathbb{N}}$. We obtain a finite set $\mathcal{N}$ of indices such that

$$\delta(\{\mathcal{K}_{B_k}\}_{k\in\mathcal{N}}) < \delta/2, \tag{11.3}$$

and

$$\delta(\{\mathcal{K}_{B_k}\}_{k\in\mathcal{N}\setminus\{n\}}) \geq \delta/2, \qquad \forall n \in \mathcal{N}. \tag{11.4}$$

By (11.3) there exist a vector $f_n \in \mathcal{K}_{B_n}$, $\|f_n\| = 1$ that can be approximated by a vector $f^n \in H^2$ of the form

$$f^n = \sum_{k\in\mathcal{N}\setminus\{n\}} f_k, \qquad f_k \in \mathcal{K}_{B_k},$$

so that

$$\|f_n - f^n\| < \delta/2.$$

Recall that $\sigma_k$ denotes the zero set of the Blaschke product $B_k$. The vectors $f_k \in \mathcal{K}_{B_k}$ can be represented in the form

$$f_k = \sum_{\lambda\in\sigma_k} c_\lambda k_\lambda, \qquad c_\lambda \in \mathbb{C}, \ k_\lambda(z) = (1-|\lambda|^2)^{1/2}(1-\overline{\lambda}z)^{-1}.$$

By the hypotheses of the lemma the system of vectors $\{k_\lambda\}_{\lambda\in\sigma_k}$ is an unconditional basis and its measure of non-orthogonality $\|R\|\cdot\|R^{-1}\|$ is at most $C$. By (11.4) the system of subspaces $\{\mathcal{K}_{B_k}\}_{k\in\mathcal{N}\setminus\{n\}}$ is a Riesz basis and its measure of non-orthogonality is at most $CV(\delta/2)$. Therefore the family of vectors $\{k_\lambda\}_{\lambda\in\sigma}$, $\sigma = \cup_{k\in\mathcal{N}\setminus\{n\}}\sigma_k$ is a Riesz basis in $\mathrm{span}(\mathcal{K}_{B_k} : k \in \mathcal{N}\setminus\{n\})$ with the measure of non-orthogonality at most $C \cdot CV(\delta/2)$.

So, if $R = R(\{k_\lambda\}_{\lambda\in\sigma})$ is a normalized orthogonalizer of $\{k_\lambda\}_{\lambda\in\sigma}$ we have

$$\|R\| \cdot \|R^{-1}\| \leq C \cdot CV(\delta/2). \tag{11.5}$$

Let $f_k = \sum_{\lambda\in\sigma_k} c_\lambda k_\lambda$ be the expansion of $f_k$ in the basis $\{k_\lambda\}_{\lambda\in\sigma_k}$. Since we assume that the family of subspaces $\{\mathcal{K}_{B_k}\}_{k\in\mathcal{N}\setminus\{n\}}$ is a Riesz basis, the zero sets $\sigma_k$ are disjoint, and we can write

$$f^n = \sum_{\lambda\in\sigma} c_\lambda k_\lambda.$$

(recall that $\sigma = \cup_{k\in\mathcal{N}\setminus\{n\}}\sigma_k$)

Let $\mathcal{K} := \mathrm{span}\{\mathcal{K}_{\Theta_k,\Delta_k,P_k} : k \in \mathcal{N}\setminus\{n\}\}$. Define

$$\tilde{f}_n := P_{K_{\Theta_n,\Delta_n,P_n}} f_n \cdot e, \qquad \tilde{f}^n := P_\mathcal{K} f^n \cdot e,$$

where $e$ is the vector from the assumptions of the main lemma. By the triangle inequality

$$\|\tilde{f}_n - \tilde{f}^n\| \leq \|\tilde{f}_n - f_n \cdot e\| + \|\tilde{f}^n - f^n \cdot e\| + \|f_n e - f^n e\|.$$

The last summand is equal to
$$\|f_n - f^n\| < \delta/2.$$
The first two terms are also easy to estimate. Let us for example estimate the second one. The subspace $\mathcal{K}$ defined above is an $S^*$-invariant subspace of $\begin{pmatrix} H^2(E) \\ L^2(E_*) \end{pmatrix}$, so it can be represented as $\mathcal{K} = \mathcal{K}_{\Theta,\Delta,P}$. Hence by Corollary 6.2

$$\|\tilde{f}^n - f^n e\| = \text{dist}\{\begin{pmatrix} f^n e \\ 0 \end{pmatrix}, \mathcal{K}_{\Theta,\Delta,P}\} = \|\mathbb{P}_+ \Theta^* f^n e\| =$$
$$= \|\mathbb{P}_+ \Theta^* \sum_{\lambda \in \sigma} c_\lambda k_\lambda e\| = \|\sum_{\lambda \in \sigma} c_\lambda k_\lambda \Theta(\lambda)^* e\|.$$

Since by the hypothesis of the lemma $\|\Theta_k(\lambda)^* e\| < 2\varepsilon$ for $\lambda \in \sigma_k$, and since by Corollary 6.3

$$\|\Theta(\lambda)^* e\| = \text{dist}\{\begin{pmatrix} k_\lambda e \\ 0 \end{pmatrix}, \mathcal{K}_{\Theta,\Delta,P}\} \leq$$
$$\leq \text{dist}\{\begin{pmatrix} k_\lambda e \\ 0 \end{pmatrix}, \mathcal{K}_{\Theta_k,\Delta_k,P_k}\} = \|\Theta_k(\lambda)^* e\|$$

one can conclude that
$$\|\Theta(\lambda)^* e\| < 2\varepsilon \qquad \forall \lambda \in \sigma.$$
Applying both inequalities from Lemma 8.2 we can get

$$\|\sum_{\lambda \in \sigma} c_\lambda k_\lambda \Theta(\lambda)^* e\|^2 \leq \|R^{-1}\|^2 \sum_{\lambda \in \sigma} |c_\lambda|^2 \|\Theta(\lambda)^* e\|^2 \leq$$
$$\leq (2\varepsilon)^2 \|R^{-1}\|^2 \sum_{\lambda \in \sigma} |c_\lambda|^2 \leq (2\varepsilon)^2 \|R^{-1}\|^2 \|R\|^2 \|\sum_{\lambda \in \sigma} c_\lambda k_\lambda\|$$

So
$$\|\tilde{f}^n - f^n e\| \leq 2\varepsilon \cdot \|R\| \cdot \|R^{-1}\| \leq 2\varepsilon \cdot C \cdot CV(\delta/2) \leq \delta/5.$$
Similarly
$$\|\tilde{f}_n - f_n e\| \leq 2\varepsilon \cdot C \leq 2\varepsilon \cdot C \cdot CV(\delta/2) \leq \delta/5.$$
(the inequality in the middle holds because the measure of non-orthogonality is always at least 1). Gathering estimates for all three summands together we have
$$\|\tilde{f}_n - \tilde{f}^n\| < \frac{9}{10}\delta,$$
$$\|\tilde{f}_n\| = (1 - \|\tilde{f}_n - f_n e\|^2)^{1/2} \geq (1 - (\delta/5)^2)^{1/2} > \sqrt{24/25} > 9/10,$$
and so
$$\left\| \|\tilde{f}_n\|^{-1} \cdot \tilde{f}_n - \|\tilde{f}_n\|^{-1} \cdot \tilde{f}^n \right\| < \delta,$$

which contradicts the definition of $\delta$ since $\tilde{f}_n \in \mathcal{K}_{\Theta_n,\Delta_n,P_n}$ and

$$\tilde{f}^n \in \mathcal{K} = \mathrm{span}\{\mathcal{K}_{\Theta_k,\Delta_k,P_k} : k \neq n\}. \qquad \square$$

## Appendix A. Proof of Bourgain's theorem on Carleson contour (Theorem 9.1).

There are three main reasons for the inclusion of this proof in the paper. First of all, to make the paper self-contained. Second, Bourgain prove this theorem only for inner functions. Although his proof works in general case, formally the theorem we need was not in [16]. And last, there was an inaccuracy in the proof, so we would like to present the corrected one.

Let $\varphi \in H^\infty$, $\|\varphi\|_\infty \leq 1$. The function $\varphi$ can be represented (see [15]) as

$$\varphi(z) = \Big(\prod_n \frac{|\lambda_n|}{\lambda_n} \cdot \frac{\lambda_n - z}{1 - \overline{\lambda}_n z}\Big) \cdot \exp\Big\{-\int_{\mathbb{T}} \frac{\xi + z}{\xi - z} d\mu(\xi)\Big\}, \tag{A.1}$$

where $\mu$ is some positive Borel measure. The exponent gives the product of outer part of $\varphi$ and its singular inner part. The absolutely continuous part of $\mu$ has density $\log|\varphi|$.

We will call the measure $\nu := \mu + (1/2) \cdot \sum_n (1 - |\lambda_n|^2)\delta_{\lambda_n}$ the representing measure of $\varphi$

The following lemma is very well known. It says that the function is small at a point $z$ if either $z$ is close to the zero of a Blaschke product part of $\varphi$ of if the Poisson potential of the representing measure $\nu$ is big.

**Lemma A.1.** *Let a function $\phi \in H^\infty$ and its representation* (A.1) *are given. Then*

$$\int_{\mathrm{clos}\,\mathbb{D}} \frac{(1-|z|^2)}{|1-\overline{\xi}z|^2} d\nu(\xi) \leq -\log|\varphi(z)| \leq 2\log\frac{1}{\varepsilon} \int_{\mathrm{clos}\,\mathbb{D}} \frac{(1-|z|^2)}{|1-\overline{\xi}z|^2} d\nu(\xi)$$

*provided* $\inf_n |b_{\lambda_n}(z)| \geq \varepsilon$.

*Proof.* Clearly

$$|\log \varphi(z)| = \int_T \frac{(1-|z|^2)}{|1-\overline{\xi}z|^2} d\mu(\xi) + \frac{1}{2}\sum_n \log\left|\frac{\lambda_n - z}{1 - \overline{\lambda}_n z}\right|^2.$$

Since $\left|\dfrac{\lambda_n - z}{1 - \overline{\lambda}_n z}\right|^2 = 1 - \dfrac{(1-|z|^2)(1-|\lambda_n|^2)}{|1-\overline{\lambda}_n z|^2}$ and

$$t \leq -\log(1-t) \leq 2\big(\log\frac{1}{\varepsilon}\big)t \qquad \text{if} \qquad 0 \leq t \leq 1-\varepsilon,$$

one can conclude

$$\int_T \frac{(1-|z|^2)}{|1-\bar{\xi}z|^2}d\mu(\xi) + \frac{1}{2}\sum \frac{(1-|z|^2)(1-|\lambda_n|^2)}{|1-\bar{\lambda}_n z|^2} \leq -\log|\varphi(z)| \leq$$
$$\leq \int_T \frac{(1-|z|^2)}{|1-\bar{\xi}z|^2}d\mu(\xi) + \log\frac{1}{\varepsilon}\sum \frac{(1-|z|^2)(1-\lambda_n|^2)}{|1-\bar{\lambda}_n z|^2}$$
$\square$

For an arc $I \subset \mathbb{T}$ let $Q(I)$ denote the Carleson square in the closed unit disk $\mathrm{clos}\,\mathbb{D}$ with base $I$, $Q(I) := \{z \in \mathrm{clos}\,\mathbb{D} : z/|z| \in I, |z| \geq 1-|I|\}$, and for $k > 0$ let $k \cdot I$ denote an ark of length $k|I|$ with the same center as $I$.

**Lemma A.2.** *Let $\varphi$ be a function in $H^\infty$, $\|\varphi\|_\infty \leq 1$, and let $\nu$ be its representing measure. Let $I$ be an arc of $\mathbb{T}$, and $\lambda$ be a point in the upper half of the Carleson square $Q(I)$ (that means $\lambda \in Q(I)$, $|\lambda| \geq 1 - |I|/2$), such that $|\varphi(\lambda)| \geq \varepsilon > 0$. Given a positive constant $M < \infty$ there exists a collection $\{I_k\}$ of disjoint open subarcs of $5I$ satisfying*

1. $\displaystyle\sum_k |I_k| \leq C_1 \big(\log\frac{1}{\varepsilon}\big) M^{-1} \cdot |I|$;
2. *If $z \in Q(I) \setminus \big(\cup_k Q(I_k)\big)$ and $\inf_n |b_{\lambda_n}(z)| \geq \gamma$ then*

$$\log|\varphi(z)|^{-1} \leq C_2\big(\log\frac{1}{\gamma}\big)\big(M + \log\frac{1}{\varepsilon}\big).$$

*Proof.* Let $\mathcal{V}$ be the collection of all open (in $5I$) subarcs $J$ of $5I$ satisfying $\nu(Q(J)) > M|I|$. By Vitali Covering Lemma there exist at most countably many disjoint intervals $J_k \in \mathcal{V}$ such that $\cup_{J \in \mathcal{V}} J \subset \cup_k 5J_k$. The open set $\cup_k 5J_k$ is a union of disjoint intervals $I_k$, which we are going to prove are the intervals $I_k$ from the conclusion of the lemma.

By Lemma A.1 we have for the point $\lambda$ from the hypothesis of the lemma

$$M\sum_k |I_k| \leq 5M\sum_k |J_k| \leq 5\sum_k \nu(Q(I_k)) \leq$$
$$\leq C|I|\int_{\cup_k Q(I_k)} \frac{(1-|z|^2)(1-|\lambda|^2)}{|1-\bar{\lambda}z|^2}d\mu(z) \leq C|I|\log|\varphi(z)|^{-1} \leq$$
$$\leq C|I|\log\frac{1}{\varepsilon},$$

so the statement 1 is proved.

Let now $z \in Q(I) \setminus \big(\cup_k Q(I_k)\big)$. Let $\tilde{I}_0$ be the interval of length $(1-|z|)$ centered at $z/|z|$. Denote $\tilde{I}_k := 2^k \tilde{I}_0$, and let $N$ be the smallest integer such that $2I \subset \tilde{I}^N$.

By the construction $\nu(Q(\tilde{I}_k)) \leq M|\tilde{I}_k|$ and taking into account that

$$\frac{1-|z|^2}{|1-\overline{z}\xi|^2} \leq C\frac{|\tilde{I}_0|}{|\tilde{I}_k|^2}, \qquad \text{for } \xi \in Q(\tilde{I}_{k+1}) \setminus Q(\tilde{I}_k)$$

one can obtain

$$\int_{Q(\tilde{I}_N)} \frac{1-|z|^2}{|1-\overline{z}\xi|^2} d\nu(\xi) \leq \int_{Q(\tilde{I}_2)} \frac{1-|z|^2}{|1-\overline{z}\xi|^2} d\nu(\xi) +$$

$$+ \sum_{k=2}^{N-1} \int_{Q(\tilde{I}_{k+1}) \setminus Q(\tilde{I}_k)} \frac{1-|z|^2}{|1-\overline{z}\xi|^2} d\nu(\xi) \leq CM + \sum_{k=2}^{N-1} CM2^{-k} \leq CM$$

If $\xi \notin Q(2I)$ then clearly

$$\frac{1-|z|^2}{|1-\overline{z}\xi|^2} \leq C\frac{1-|\lambda|^2}{|1-\overline{\lambda}\xi|^2},$$

where $C$ is an absolute constant. Hence by Lemma A.1 and hypothesis

$$\int_{\operatorname{clos}\mathbb{D}\setminus Q(\tilde{I}_N)} \frac{1-|z|^2}{|1-\overline{z}\xi|^2} d\nu(\xi) \leq$$

$$\leq C \int_{\operatorname{clos}\mathbb{D}\setminus Q(\tilde{I}_N)} \frac{1-|\lambda|^2}{|1-\overline{\lambda}\xi|^2} d\nu(\xi) \leq C\log|\varphi(\lambda)|^{-1} \leq C\log\frac{1}{\varepsilon}$$

$\square$

**Lemma A.3.** *In the situation of Lemma A.2 and given $\gamma > 0$ there exists a set $\mathcal{O}(I)$ in $Q(I)$ with rectifiable boundary $\Gamma(I)$ such that*
1. *The arc-length on $\Gamma(I)$ is a Carleson measure with Carleson norm at most $C_3\gamma[M + \log(1/\varepsilon)]$*
2. $|\varphi(z)| < \gamma$ *if* $z \in \mathcal{O}(I)$;
3. $|\varphi(z)| > \exp\left\{-C_2\left(\log\frac{1}{\gamma}\right)\left(M + \log\frac{1}{\varepsilon}\right)\right\}$ *if* $z \in Q(I) \setminus \{\mathcal{O}(I) \cup \bigcup_k Q(I_k)\}$

*Proof.* Define

$$\mathcal{O}(I) := \left[Q(I) \setminus \bigcup_k Q(I_k)\right] \cap \bigcup_{\lambda_n \in Q(2I)} \{z \in \mathbb{D} : |b_{\lambda_n}(z)| < \gamma\},$$

where $\{\lambda_n\}$ are zeroes of $\varphi$. Then the statement 2 of the lemma is trivial. Statement 3 follows from Statement 2 of Lemma A.2.

It is easy to see that

$$|\Gamma| \leq C\gamma \sum_{\lambda_n \in Q(2I)} (1 - |\lambda_n|^2) \leq C\gamma \nu(Q(2I)) \leq$$

$$\leq C\gamma |I| \int_{Q(2I)} \frac{1-|\lambda|^2}{|1-\lambda z|^2} d\nu(z) \leq C\gamma \Big(\log \frac{1}{\varepsilon}\Big)|I|$$

Suppose now $J$ is a subinterval of $5I$. If $J$ is contained in $\bigcup_k I_k$, then $\Gamma \cap Q(J) = \varnothing$. Otherwise by the construction of the intervals $I_k$

$$\sum_{\lambda_n \in Q(2J)} (1 - |\lambda_n|^2) \leq \nu(2J) \leq 2M|J|.$$

Therefore

$$|\Gamma \cap Q(J)| \leq C\gamma \sum_{\alpha_n \in Q(2J)} (1 - |\lambda_n|^2) \leq C\gamma M|J|$$

□

*Proof of Theorem 9.1.* The proof of Theorem 9.1 follows now a standard reasoning of constructing generations of intervals in $\mathbb{T}$. Take $M = 100C_1 \log \frac{1}{\varepsilon}$, $\gamma = \min\{\varepsilon, 1/(2C_3[M + \log(1/\varepsilon)])\}$ and

$$\varepsilon'(\varepsilon) = \exp\Big(-C_2 \log \frac{1}{\gamma}\Big(M + \log \frac{1}{\varepsilon}\Big)\Big)$$

Start with first generation $\mathcal{F}^{(1)} = \{\mathbb{T}\}$. Assume now the generation $\mathcal{F}^{(s)}$ is constructed and $J$ be an interval in $\mathcal{F}^{(s)}$. Denote by $\mathcal{D}(J)$ the set of all dyadic subintervals $I$ of $J$ such that the upper half of the Carleson square $Q(I)$ contains a point $\lambda$ where $|\varphi(\lambda)| \geq \varepsilon$. Therefore

$$|\varphi(z)| \leq \varepsilon \qquad \text{on} \qquad Q(J) \setminus \bigcup_{I \in \mathcal{D}(J)} Q(I). \tag{A.2}$$

To each $I \in \mathcal{D}(J)$ we apply lemmas A.2 and A.3 providing intervals $\{I_k\}$ and a region $\mathcal{O}(I)$ of $Q(I)$. Define

$$\mathcal{F}^{(s+1)} := \bigcup_{J \in \mathcal{F}^{(s)}} \bigcup_{I \in \mathcal{D}(J)} \{I_k : k = 1, 2, \ldots\}.$$

Note that intervals in $\mathcal{F}^{(s+1)}$ are not necessary disjoint.

By the construction for each $I \in \mathcal{D}(J)$

$$\sum_k |I_k| \leq \frac{1}{100}|I| \tag{A.3}$$

For $J \in \mathcal{F}^{(s)}$ let $R_s(J)$ be the union of $Q(J) \setminus \bigcup_{I \in \mathcal{D}(J)} Q(I)$ and $\bigcup_{I \in \mathcal{D}(J)} \mathcal{O}(I)$. Take $R_s := \bigcup_{J \in \mathcal{F}^{(s)}} R_s(J)$. Thus by (A.2) and statement 2 of Lemma A.3 imply $|\varphi(z)| \leq \varepsilon$ on $R_s(J)$ and therefore on $R_s$.

Also by statement 3 of Lemma A.3

$$|\varphi(z)| \geq \varepsilon'(\varepsilon) \quad \text{if} \quad z \in \bigcup_{J \in \mathcal{F}^{(s)}} R(J) \setminus \left( R_s \cup \bigcup_{J \in \mathcal{F}^{(s+1)}} R(J) \right)$$

Define $\mathcal{O} = \bigcup R_s$. Our construction yields $|\varphi| \leq \varepsilon$ on $\mathcal{O}$ and $|\varphi| \geq \varepsilon'(\varepsilon)$ on $\mathbb{D} \setminus \mathcal{O}$.

The reader can easily see that (A.3) and statement 1 of Lemma A.3 imply that the arc-length on $\partial \mathcal{O}$ is a Carleson measure with Carleson norm at most 10. □

## References

[1] N. K. Nikolskii, *Treatise on the shift operator*, Springer-Verlag, NY etc. 1985.
[2] V. E. Katsnelson, *On basis conditions for a system of principal vectors of certain classes of operators*, Funct. Anal. Appl., **1**(1967).
[3] I. C. Gohberg and M. G. Krein, *Introduction to the theory of linear nonselfadjoint operators*, Providence, R.I., American Mathematical Society, 1969.
[4] B. Sz.-Nagy and C. Foiaş, *Harmonic analysis of operators in Hilbert space*, Akadémia Kiadó, Budapest, 1970.
[5] R. Hunt, B. Muckenhoupt and R. Wheeden, *Weighted norm inequalities for conjugate function and Hilbert transform*, Trans. Amer. Math. Soc., **176** (1973), 227–251.
[6] N. K. Nikolskii and B. S. Pavlov, *Expansions in eigenfunctions of non-unitary operators and the characteristic function*, Zap. Nauchn. Semin. LOMI, **11**(1968), 150–203 (Russian).
[7] S. R. Treil *A spatially compact system of eigenvectors forms a Riesz basis if it is uniformly minimal*, Dokl. Akad. Nauk SSSR, **288** (1986), 308 – 312; (Russian) English translation: Soviet Math. Dokl., v.33 (1986), 675 – 679.
[8] S. R. Treil, *Geometric methods in spectral theory of vector valued functions: Some recent results*, Operator Theory: Adv. Appl., **42** (1989), pp. 209 – 280.
[9] S. R. Treil, *Hankel operators, imbedding theorems, and bases of the invariant subspaces of the multiple shift*, Algebra i Analiz, **1** (1989), No 6, 200–234 (Russian); english translation: St. Petersburg Math. Journal, **1**(1990), 1515–1548.
[10] V. I. Vasyunin, *Unconditionally convergent spectral decompositions and interpolation problems*, Trudy Matem. Inst. Akad. Nauk SSSR, **130** (1978), 5–49 (Russian); English translation: Proc. Steklov Inst. Math., Issue 4 (1979), 1–53.
[11] C. Caratheodory, *Theory of functions of complex variable*, New York, Chelsea Publ. Co., 1954.
[12] I. I. Privalov, *Introduction to the theory of analytic functions*, GIFML, Moskow, 1960.
[13] F. F. Bonsall, *Boundedness of Hankel matrices*, J. London math. Soc., **29** (1984), Part 2, 289—300.
[14] S. V. Khruščev and N. K. Nikolskii, *Functional model and some problems of spectral theory of functions*, Trudy Mat. Inst. Steklov, **176** (1987), 97–209 (Russian); English traslation: Proc. Steklov Inst. Math.
[15] J. B. Garnett, *Bounded analytic functions*, Academic Press, New York, 1981.
[16] J. Bourgain, *On finitely generated closed ideals in $H^\infty(\mathbb{D})$*, An. Institute Fourier, Grenoble, **35**(1985), No 4, 163–174.

Serguei Treil, Department of Mathematics, Michigan State University, East Lansing, MI 48823

*E-mail address*: `treil@@math.msu.edu`